\documentclass[11pt,a4paper]{article}

\usepackage[top=2.5cm,bottom=2.5cm,right=2.5cm,left=2.5cm]{geometry}
\usepackage[T1]{fontenc}
\usepackage{lmodern}
\usepackage{amsmath,amsthm}
\usepackage{amssymb,amsfonts}
\usepackage{xcolor}
\usepackage{mathtools,bm}
\usepackage{amstext}
\usepackage{mathrsfs}
\usepackage{multirow}
\usepackage{array}
\usepackage{algpseudocode}
\usepackage{amsopn}
\usepackage{threeparttable}
\usepackage{enumitem}
\usepackage{subfigure}
\usepackage{cite}

\usepackage{rotating}
\usepackage{fancyvrb}
\usepackage{graphicx}
\usepackage{algorithm}
\usepackage{ifpdf}
\usepackage[hidelinks]{hyperref}
\newtheorem{theorem}{Theorem}
\newtheorem{lemma}[theorem]{Lemma}

\newtheorem{definition}[theorem]{Definition}
\newtheorem{assumption}[theorem]{Assumption}
\newenvironment{keywords}{\par\begingroup\small\noindent\textbf{Keywords.}\ }{\par\endgroup\smallskip}
\newenvironment{MSCcodes}{\par\begingroup\small\noindent\textbf{MSC codes.}\ }{\par\endgroup\smallskip}
\newcommand{\email}[1]{\href{mailto:#1}{\nolinkurl{#1}}}

\newcommand{\dist}{\operatorname{dist}}

\providecommand{\argmin}{\operatorname*{arg\,min}}

\newcommand{\TheTitle}{ALDER: A Single-Loop Augmented-Lagrangian Damped Extragradient Method for Nonconvex--Concave Minimax Problems with Nonlinear Coupled Constraints}
\newcommand{\TheAuthors}{H. Zhang, M. Zhang, Z. Xu, and Y.-H. Dai}

\title{\TheTitle\thanks{This work was supported by the National Key R\&D Program of China (Grant Nos.~2025YFA1017801 and 2025YFA1017800) and the National Natural Science Foundation of China (Grant No.~12471294).}}

\author{
	Huiling Zhang\thanks{Department of Mathematics, Shanghai University, Shanghai 200444, People’s Republic of China; LSEC, ICMSEC, Academy of Mathematics and Systems Science, Chinese Academy of Sciences, Beijing 100190, China (\email{zhanghl1209@shu.edu.cn}).}
	\and
	Minhao Zhang\thanks{Department of Mathematics, Shanghai University, Shanghai 200444, People’s Republic of China (\email{zhangminhao@shu.edu.cn}).}
	\and
	Zi Xu\thanks{Department of Mathematics, Shanghai University, Shanghai 200444, People’s Republic of China; Newtouch Center for Mathematics of Shanghai University, Shanghai 200444, People’s Republic of China (\email{xuzi@shu.edu.cn}). Corresponding author.}
	\and
	Yu-Hong Dai\thanks{State Key Laboratory of Mathematical Sciences, Academy of Mathematics and Systems Science, Chinese Academy of Sciences, Beijing 100190, China; School of Mathematical Sciences, University of Chinese Academy of Sciences, Beijing 100049, China. (\email{dyh@lsec.cc.ac.cn}).}
}
\date{}

\ifpdf
\hypersetup{
	pdftitle={\TheTitle},
	pdfauthor={\TheAuthors}
}
\fi

	\let\siamproof\proof
	\let\endsiamproof\endproof
	\renewenvironment{proof}[1]{\siamproof}{\endsiamproof}
	
\begin{document}
		
		\maketitle
		
\begin{abstract}
	ALDER, a single-loop augmented-Lagrangian damped extragradient
	method, is proposed for smooth nonconvex--concave minimax problems
	with nonlinear coupled inequality constraints.
	Within an augmented-Lagrangian framework, ALDER combines
	prediction--correction steps for the primal variables with
	explicit projected-and-relaxed multiplier updates in a single loop.
	The penalty parameter is fixed throughout each run, eliminating
	inner-loop subproblem solves and iterative penalty updates.
	We establish uniform boundedness of the relaxed multiplier sequence
	and, under an extended Mangasarian--Fromovitz constraint qualification,
	boundedness of the projected output multiplier.
	Building on these bounds, we show that ALDER finds an
	$\varepsilon$-KKT point within
	$\mathcal{O}(\varepsilon^{-7/2})$ complete first-order oracle
	evaluations in the strongly concave setting and
	$\mathcal{O}(\varepsilon^{-4})$ in the concave setting. These bounds improve upon the best known complexity bounds
	for multi-loop algorithms for this problem class, achieved by FAL,
	by removing a $\log(1/\varepsilon)$ factor, while retaining
	a single-loop update structure.
\end{abstract}

		
		
		\section{Introduction}\label{sec1}
		Consider the following smooth nonconvex minimax optimization problem with coupled nonlinear inequality constraints:
		\begin{equation}\label{prob:P}
			\min_{x\in\mathcal{X}}\max\limits_{\substack{ y\in\mathcal{Y}\\c(x,y)\le0}}f(x,y),\tag{P}
		\end{equation}
		\addtocounter{equation}{1}
		where $\mathcal{X}\subseteq \mathbb{R}^n$ and $\mathcal{Y}\subseteq \mathbb{R}^m$ are nonempty convex and compact sets, $f(x,y): \mathcal{X} \times \mathcal{Y} \rightarrow \mathbb{R}$ is a smooth function and possibly nonconvex in $x$, and concave in $y$. The mapping $c=(c_1,\ldots,c_p)^\top$ is smooth on a neighborhood of $\mathcal X\times\mathcal Y$, and each $c_i(x,\cdot)$ is convex. We use $c_i(x,y)\le0$ to denote the $i$th constraint, and define $\mathcal{K}:= \{1,\dots , p\}$ as the index set of the constraints.
		
		When the coupled inequalities are absent, \eqref{prob:P} degenerates into the following special minimax optimization problem:
		\begin{equation}\label{prob:uncoupled}
			\min \limits_{x\in \mathcal{X}}\max \limits_{y\in \mathcal{Y}}f(x,y),
		\end{equation}
		which has been a focus of research in machine learning, signal processing, statistics, economics and many other research fields \cite{Berger,Facchinei,Ghosh,Sanjabi,Wu,Yang} recently. Many applications can be formulated into \eqref{prob:uncoupled}, including the power control and transceiver design problem in signal processing \cite{Lu}, reinforcement  learning \cite{Cai,Moriarty,Qiu}, the adversarial attack to deep neural networks (DNNs) \cite{Chen}, distributed nonconvex optimization \cite{Giannakis,Mateos}, robust learning over multiple domains \cite{Qian}.  We call \eqref{prob:uncoupled} a (non)convex-(non)concave minimax problem if $f(x, y)$ is (non)convex in $x$ and (non)concave in $y$. Minimax problems have been extensively studied in the (strongly) convex--(strongly) concave \cite{Lan2016,Mokhtari,Nemirovski,Nesterov2007,Thek2022}, nonconvex--(strongly) concave \cite{Lin2019,Lin2020,Lu,Pan,Xu,Xu21,Zhang}, and nonconvex--nonconcave \cite{Bohm,DiakonikolasEG,Doan,Lee,Sanjabi18,SongOptDE} settings.
		Zhang et al. \cite{Zhang2021} developed a generic acceleration framework for the nonconvex--strongly concave minimax problem \eqref{prob:uncoupled}, attaining the optimal complexity bound of $\mathcal{O}(\sqrt{\kappa_y}\varepsilon^{-2})$ for finding an $\varepsilon$-stationary point \cite{Li21}.
		Zhang and Xu \cite{Zhang2026} proposed a single-loop projected damped extragradient method for smooth nonconvex--concave minimax optimization with a gradient complexity of $\mathcal{O}(\varepsilon^{-5/2})$, the best known among single-loop first-order methods for this setting.
		
		Compared with \eqref{prob:uncoupled}, problem \eqref{prob:P} models a broader range of applications, including adversarial attacks in resource allocation and network flow problems \cite{Tsaknakis}. Finding a global solution to \eqref{prob:P} is NP-hard even when $f(x,y)$ is strongly convex in $x$ and strongly concave in $y$ \cite{Tsaknakis}.
		For smooth minimax problems with coupled linear inequality constraints, Tsaknakis et al. \cite{Tsaknakis} proposed a family of multiplier gradient descent (MGD) algorithms with a total first-order complexity of $\mathcal{O}(\varepsilon^{-2}\log(1/\varepsilon))$ in the strongly convex--strongly concave setting.
		Dai et al. \cite{Dai} developed a proximal gradient multi-step ascent descent (PGmsAD) method for nonsmooth convex--concave minimax problems with bilinear coupling and linear equality constraints, achieving a total first-order complexity of $\mathcal{O}(\varepsilon^{-2}\log(1/\varepsilon))$.  For nonconvex--concave minimax problems with nonlinear coupled inequality constraints, Lu and Mei \cite{lu23} proposed a first-order augmented Lagrangian (FAL) method that finds an $\varepsilon$-KKT point with a total first-order complexity of $\mathcal{O}(\varepsilon^{-4}\log(1/\varepsilon))$. They subsequently obtained a total first-order complexity bound of $\mathcal{O}(\varepsilon^{-7/2}\log(1/\varepsilon))$ by exploiting strong concavity in the inner variable \cite{LuMei2026SC}. Hu et al. \cite{HuTohWangXiao2024} studied nonconvex--strongly concave minimax problems with nonlinear conic constraints and proposed an MMPen-based single-loop subgradient descent--ascent method, without an explicit $\varepsilon$-KKT oracle complexity bound for this method. Cao et al. \cite{CaoZengZhangZhou2026} proposed SPACO, a single-loop stochastic penalty method for nonconvex--concave minimax problems with nonlinear inequality constraints. They established an iteration bound of $\mathcal{O}(\varepsilon^{-1/\tau})$ for an expected reduced KKT residual that excludes complementarity, where $\tau$ depends on the stepsize and penalty exponents. Table~\ref{tab:coupled-comparison} summarizes representative methods for problems with coupled constraints. 
		\begin{table}[!ht]
			\centering
			\caption{Complexity comparison for minimax problems with coupled constraints.}
			\label{tab:coupled-comparison}
			\begingroup
			\small
			\setlength{\tabcolsep}{3pt}
			\renewcommand{\arraystretch}{1.15}
			\begin{tabular}{@{}
					>{\raggedright\arraybackslash}p{.18\linewidth}
					>{\raggedright\arraybackslash}p{.10\linewidth}
					>{\raggedright\arraybackslash}p{.17\linewidth}
					>{\raggedright\arraybackslash}p{.08\linewidth}
					>{\raggedright\arraybackslash}p{.20\linewidth}
					>{\raggedright\arraybackslash}p{.16\linewidth}@{}}
				\hline
				Method & Setting & Constraints & Loops & Complexity & Criterion \\
				\hline
				MGD \cite{Tsaknakis}
				& SC--SC & Linear ineq. & Multi
				& $\mathcal{O}(\varepsilon^{-2}\log(1/\varepsilon))$
				& $\varepsilon$-stationarity \\
				
				PGmsAD \cite{Dai}
				& C--C\textsuperscript{*} & Linear eq. & Multi
				& $\mathcal{O}(\varepsilon^{-2}\log(1/\varepsilon))$
				& $\varepsilon$-KKT \\
				\hline
				FAL \cite{LuMei2026SC}
				& NC--SC & Nonlinear ineq. & Multi
				& $\mathcal{O}(\varepsilon^{-7/2}\log(1/\varepsilon))$
				& $\varepsilon$-KKT \\
				
				\textbf{ALDER}
				& NC--SC & Nonlinear ineq. & \textbf{Single}
				& $\boldsymbol{\mathcal{O}(\varepsilon^{-7/2})}$
				& $\varepsilon$-KKT \\
				\hline
				FAL \cite{lu23}
				& NC--C & Nonlinear ineq. & Multi
				& $\mathcal{O}(\varepsilon^{-4}\log(1/\varepsilon))$
				& $\varepsilon$-KKT \\
				
				\textbf{ALDER}
				& NC--C & Nonlinear ineq. & \textbf{Single}
				& $\boldsymbol{\mathcal{O}(\varepsilon^{-4})}$
				& $\varepsilon$-KKT \\
				\hline
			\end{tabular}
			
			\smallskip
			\begin{minipage}{\linewidth}
				\footnotesize
				NC: nonconvex; C: convex/concave; SC: strongly convex/strongly concave,
				according to the variable.
				\textsuperscript{*}Bilinear coupling.
			\end{minipage}
			\endgroup
		\end{table}

		\textbf{Contributions.}
		We develop a single-loop first-order method for smooth nonconvex--concave minimax problems with nonlinear coupled inequality constraints. The main contributions are as follows.
		\begin{itemize}
			\item We propose ALDER, a single-loop augmented-Lagrangian damped extragradient method with relaxed multiplier updates. The method combines prediction--correction steps for both primal variables with an explicit projected-and-relaxed multiplier update, avoiding inner-loop subproblem solves and penalty-parameter updates.
			
			\item We establish uniform boundedness of the relaxed multiplier sequence and, under an extended Mangasarian--Fromovitz constraint qualification, boundedness of the projected output multiplier. These bounds provide the multiplier control needed to establish $\varepsilon$-KKT guarantees, including stationarity in both primal variables, feasibility, and componentwise approximate  complementarity, with exact dual nonnegativity.
			
			\item Under the corresponding assumptions, ALDER finds an $\varepsilon$-KKT point using $\mathcal{O}(\varepsilon^{-7/2})$ and $\mathcal{O}(\varepsilon^{-4})$ complete first-order evaluations in the nonconvex--strongly concave and nonconvex--concave settings, respectively. To the best of our knowledge, these bounds have the best dependence on $\varepsilon$ among existing first-order methods for this problem class, including both single-loop and multi-loop methods. In particular, ALDER removes the logarithmic factor in the corresponding multi-loop FAL bounds of $\mathcal{O}(\varepsilon^{-7/2}\log(1/\varepsilon))$ \cite{LuMei2026SC} and $\mathcal{O}(\varepsilon^{-4}\log(1/\varepsilon))$ \cite{lu23}.
		\end{itemize}
		
		\textbf{Organization.}
		Section~\ref{sec:method} presents ALDER{\color{blue}.}
		Section~\ref{sec:analysis}  establishes the complexity bounds in the strongly concave and concave settings.
		Section~\ref{sec:nu} presents numerical experiments illustrating the performance of ALDER.
		Section~\ref{section_conclu} concludes the paper.
		
		\textbf{Notation.}
		We use $\|\cdot\|$ for the Euclidean norm and the induced matrix norm, and $\langle\cdot,\cdot\rangle$ for the Euclidean inner product.
		For a vector $x$, $[x]_i$ denotes its $i$th component, and $(x;y)$ denotes the vertical concatenation of $x$ and $y$.
		We write $\mathbb R^p$ for the $p$-dimensional Euclidean space and $\mathbb R_+^p$ for its nonnegative orthant.
		For a differentiable function $f(x,y,z)$, $\nabla_x f$, $\nabla_y f$, and $\nabla_z f$ denote its partial gradients, and $\nabla_{x,z}f:=(\nabla_x f;\nabla_z f)$.
		For a differentiable constraint mapping $c(x,y)$, its Jacobian is $Jc=[J_xc\ J_yc]$, where $J_xc$ and $J_yc$ are the partial Jacobians.
		
		For a nonempty closed convex set $\mathcal C$, $\mathcal P_{\mathcal C}$ denotes the Euclidean projection onto $\mathcal C$, and its normal cone at $u\in\mathcal C$ is $N_{\mathcal C}(u):=\{v:\langle v,w-u\rangle\le 0 \text{ for all } w\in\mathcal C\}$.
		We write $\dist(u,\mathcal C):=\inf_{v\in\mathcal C}\|u-v\|$ and $[a]_+:=\mathcal P_{\mathbb R_+^p}(a)$ for $a\in\mathbb R^p$.
		
		A continuously differentiable function $f$ is $L$-smooth on $\mathcal X$ if $\|\nabla f(x)-\nabla f(y)\|\le L\|x-y\|$ for all $x,y\in\mathcal X$ and some $L>0$.
		It is $\mu$-strongly concave on $\mathcal X$ if, for some $\mu>0$,
		\[
		f(y)\le f(x)+\langle\nabla f(x),y-x\rangle
		-\frac{\mu}{2}\|y-x\|^2,
		\qquad x,y\in\mathcal X.
		\]
		The notation $\mathcal O(\cdot)$ suppresses constants independent of the target accuracy, and $\tilde{\mathcal O}(\cdot)$ additionally suppresses logarithmic factors.
		
		\section{An Augmented Lagrangian Damped Extragradient Method}
		\label{sec:method}
		
		We propose ALDER, a single-loop augmented Lagrangian damped extragradient method with relaxed multiplier updates for problem \eqref{prob:P}.
		The method combines projected prediction--correction steps for both primal variables with explicit multiplier and auxiliary-state updates.
		
		The augmented Lagrangian associated with \eqref{prob:P} is
		\begin{equation}\label{def:qgamma}
			\mathcal L_\gamma(x,y,\lambda)
			:=f(x,y)-P_\gamma(x,y,\lambda),
		\end{equation}
		where
		\begin{equation}\label{def:Pgamma}
			P_\gamma(x,y,\lambda)
			:=\frac{\|[\lambda+\gamma c(x,y)]_+\|^2-\|\lambda\|^2}{2\gamma},
		\end{equation}
		and $\gamma>0$ is the penalty parameter.
		To regularize the primal updates, we introduce the surrogate
		\begin{equation}\label{def:G}
			G(x,y,z,\lambda)
			:=\mathcal L_\gamma(x,y,\lambda)
			+\frac{r_\gamma}{2}\|x-z\|^2
			-\frac{\rho}{2}\|y\|^2.
		\end{equation}
		The quadratic term in $x$ regularizes the descent step around a moving center $z$, while the quadratic term in $y$ provides additional concavity when $\rho>0$.
		
ALDER builds on the projected damped extragradient mechanism
of~\cite{Zhang2026} for the primal prediction--correction steps.
The method in~\cite{Zhang2026} addresses problems without coupled
constraints and therefore does not involve Lagrange multipliers.
Here, the prediction--correction steps are applied to the
augmented-Lagrangian surrogate and coupled with explicit
projected-and-relaxed multiplier updates.
Each iteration consists of the following updates.
\begin{itemize}
	\item A projected descent--ascent step generates the prediction
	point $(\widetilde x_t,\widetilde y_t)$ using the stored
	projection residuals $\xi_t,n_t$ and damping state $v_t$.
	Gradients evaluated at this point are then used to compute
	the corrected primal iterates $(x_{t+1},y_{t+1})$.
	The $y$-correction also uses
	\[
	\bar v_t
	=\alpha v_t+\vartheta
	\nabla_yG(\widetilde x_t,\widetilde y_t,z_t,\lambda_t),
	\]
	which combines the previous damping state with the gradient
	at the prediction point.
	
	\item The correction projections determine the updated
	residuals $\xi_{t+1}$ and $n_{t+1}$.
	These residuals, together with the gradient at the corrected
	point, are used to update the damping state $v_{t+1}$.
	
	\item The corrected primal point determines the projected
	multiplier
	\[
	\widetilde\lambda_{t+1}
	=[\lambda_t+\gamma c(x_{t+1},y_{t+1})]_+.
	\]
	The multiplier used in the next iteration is obtained by
	relaxation,
	\[
	\lambda_{t+1}
	=(1-\theta_t/\gamma)\lambda_t
	+(\theta_t/\gamma)\widetilde\lambda_{t+1},
	\]
	and the center is updated as
	$z_{t+1}=z_t+\beta(x_{t+1}-z_t)$.
\end{itemize}

The nonlinear coupled constraints introduce an additional
interaction between the primal steps and the multiplier updates:
the primal iterates determine the constraint residuals driving
the multiplier updates, while the multipliers enter the surrogate
gradients used in subsequent primal steps.
Controlling this interaction requires multiplier bounds and
feasibility and complementarity estimates that do not arise
in the setting of~\cite{Zhang2026}.
In particular, the relaxation separates the multiplier sequence
used to generate the iterates from the projected multiplier
used in the approximate KKT output.
The analysis in Section~\ref{sec:analysis} establishes uniform
boundedness of the former and, under an extended
Mangasarian--Fromovitz constraint qualification, boundedness
of the latter, providing the multiplier control needed for
the KKT guarantees.
Thus, the connection to~\cite{Zhang2026} lies in the primal
prediction--correction mechanism, whereas the treatment of
coupled constraints requires additional algorithmic components
and a corresponding multiplier analysis.
Unlike the multi-loop FAL methods~\cite{lu23,LuMei2026SC},
ALDER updates the primal variables, multiplier, and auxiliary
states within a single loop, with no inner subproblem solves
and with the penalty parameter fixed throughout each run.

		Algorithm~\ref{alg:cpdeg} gives the complete procedure.
		Its output rule selects an iterate using the primal stationarity residual $S_t$ and returns the associated projected multiplier.

\begin{algorithm}[!ht]
	\caption{ALDER: Augmented Lagrangian Damped Extragradient
		with Relaxed Multiplier Updates}
	\label{alg:cpdeg}
	\begin{algorithmic}
		\State \textbf{Initialization:}
		Choose $(x_0,y_0)\in\mathcal X\times\mathcal Y$ satisfying
		$c(x_0,y_0)\le0$.
		Input $\gamma>0$, $\rho\ge0$, $r_\gamma\ge0$, $h>0$, $\alpha\in(0,1)$, $\vartheta\ge0$, $\beta\in(0,1]$, an integer $T\ge1$, and $0<\theta_t\le\gamma$;
		set $\lambda_0=0$, $z_0=x_0$, and $\xi_0=n_0=v_0=0$.
		
		\For{$t=0,1,\ldots,T-1$}
		\State \textit{Prediction and damping:}
		\State
		$\widetilde x_t
		=\mathcal P_{\mathcal X}\!\left(
		x_t-h\bigl(
		\nabla_xG(x_t,y_t,z_t,\lambda_t)+\xi_t
		\bigr)\right)$.
		\State
		$\widetilde y_t
		=\mathcal P_{\mathcal Y}\!\left(
		y_t+h\bigl(
		\nabla_yG(x_t,y_t,z_t,\lambda_t)-n_t+v_t
		\bigr)\right)$.
		\State
		$\bar v_t
		=\alpha v_t+\vartheta
		\nabla_yG
		(\widetilde x_t,\widetilde y_t,z_t,\lambda_t)$.
		
		\State \textit{Primal correction:}
		\State
		$x_{t+1}
		=\mathcal P_{\mathcal X}\!\left(
		x_t-h\nabla_xG
		(\widetilde x_t,\widetilde y_t,z_t,\lambda_t)
		\right)$.
		\State
		$y_{t+1}
		=\mathcal P_{\mathcal Y}\!\left(
		y_t+h\bigl(
		\nabla_yG
		(\widetilde x_t,\widetilde y_t,z_t,\lambda_t)+\bar v_t
		\bigr)\right)$.
		
		\State \textit{Residual and damping updates:}
		\State
		$\displaystyle
		\xi_{t+1}
		=\frac{x_t-x_{t+1}}{h}
		-\nabla_xG
		(\widetilde x_t,\widetilde y_t,z_t,\lambda_t)$.
		\State
		$\displaystyle
		n_{t+1}
		=\frac{
			(y_t-y_{t+1})/h
			+\nabla_yG
			(\widetilde x_t,\widetilde y_t,z_t,\lambda_t)
			+\bar v_t
		}{1+\vartheta}$.
		\State
		$v_{t+1}
		=\alpha v_t+\vartheta\bigl(
		\nabla_yG
		(x_{t+1},y_{t+1},z_t,\lambda_t)-n_{t+1}
		\bigr)$.
		
		\State \textit{Multiplier and center updates:}
		\State
		$\widetilde\lambda_{t+1}
		=\mathcal P_{\mathbb R_+^p}\!\left(
		\lambda_t+\gamma c(x_{t+1},y_{t+1})
		\right)$.
		\State
		$\displaystyle
		\lambda_{t+1}
		=\left(1-\frac{\theta_t}{\gamma}\right)\lambda_t
		+\frac{\theta_t}{\gamma}\widetilde\lambda_{t+1}$.
		\State
		$z_{t+1}=z_t+\beta(x_{t+1}-z_t)$.
		\EndFor
		
		\State \textbf{Output:}
		Define
		$S_t:=
		\|\nabla_x\mathcal L_\gamma
		(x_{t+1},y_{t+1},\lambda_t)+\xi_{t+1}\|^2
		+\|-\nabla_y\mathcal L_\gamma
		(x_{t+1},y_{t+1},\lambda_t)
		+\rho y_{t+1}+n_{t+1}\|^2$.
		Choose $j\in\argmin_{0\le t<T}S_t$ and return
		$(x_{\rm out},y_{\rm out},\lambda_{\rm out})
		=(x_{j+1},y_{j+1},\widetilde\lambda_{j+1})$.
	\end{algorithmic}
\end{algorithm}

		\section{Complexity Analysis}\label{sec:analysis}
	In this section, we establish the iteration complexity of Algorithm~\ref{alg:cpdeg} in the nonconvex--concave and nonconvex--strongly concave settings. We first state the assumptions and define the approximate KKT criterion.
	
	\begin{assumption}\label{ass:basic}
		The following conditions hold.
		\begin{enumerate}[label=(A\arabic*)]
			\item $\mathcal X\subseteq\mathbb R^n$ and
			$\mathcal Y\subseteq\mathbb R^m$ are nonempty compact convex sets.
			Let $D_Y:=\max_{y\in\mathcal Y}\|y\|$.
			
			\item The function $f$ is continuously differentiable on a neighborhood of
			$\mathcal X\times\mathcal Y$, its full gradient is $L_f$-Lipschitz,
			and $f(x,\cdot)$ is concave for every $x$.
			Define
			\[
			G_x:=\max_{\mathcal X\times\mathcal Y}\|\nabla_xf(x,y)\|,
			\qquad
			G_y:=\max_{\mathcal X\times\mathcal Y}\|\nabla_yf(x,y)\|.
			\]
			
			\item The mapping $c$ is continuously differentiable on a neighborhood of
			$\mathcal X\times\mathcal Y$, and each $c_i(x,\cdot)$ is convex.
			There exist $C_c,J_c,L_c,L_{Jc}>0$ such that, for all
			$(x,y),(x',y')\in\mathcal X\times\mathcal Y$,
			\begin{align*}
				\|c(x,y)\|&\le C_c,
				\qquad \|Jc(x,y)\|\le J_c,\\
				\|c(x,y)-c(x',y')\|
				&\le L_c\|(x-x',y-y')\|,\\
				\|Jc(x,y)-Jc(x',y')\|
				&\le L_{Jc}\|(x-x',y-y')\|.
			\end{align*}
			
			\item There exist $\delta>0$ and $Z>0$ such that, for every
			$(x,y)\in\mathcal X\times\mathcal Y$, there is
			$d(x,y)\in N_{\mathcal Y}^*(y)$ with $\|d(x,y)\|\le Z$ satisfying
			\[
			\delta+\nabla_yc_i(x,y)^\top d(x,y)\le0
			\quad\text{for every }i\text{ with }c_i(x,y)>0,
			\]
			where $N_{\mathcal Y}^*(y):=
			\{d:\langle d,u\rangle\le0,\ \forall u\in N_{\mathcal Y}(y)\}$.
		\end{enumerate}
	\end{assumption}
	
	Conditions (A1)--(A3) specify the domain, smoothness, and concavity structure.
	Condition (A4) is an extended MFCQ condition, also used in nonconvex constrained optimization \cite{Shi}, and provides the constraint qualification needed for the output multiplier bound.
	Under (A1) and (A3), this condition is weaker than
	the uniform Slater condition imposed on the inner coupled constraints in
	\cite[Assumption~5(ii)]{lu23}.
	The strongly concave result additionally requires the following assumption.
	
	\begin{assumption}\label{ass:strong-y}
		There is a fixed $\sigma_y>0$ such that $f(x,\cdot)$ is
		$\sigma_y$-strongly concave on $\mathcal Y$ for every $x\in\mathcal X$.
	\end{assumption}

	Following~\cite[Section~3.1]{lu23}, we relate the KKT criterion
	to problem~\eqref{prob:P}. For each fixed $x\in\mathcal X$
	with a nonempty inner feasible set, the minimax theorem,
	together with compactness and convexity of $\mathcal Y$,
	concavity of $f(x,\cdot)$, and convexity of $c_i(x,\cdot)$, gives
	\begin{equation}\label{pro:kkt}
		\max_{\substack{y\in\mathcal Y\\c(x,y)\le0}}f(x,y)
		=
		\inf_{\lambda\in\mathbb R_+^p}\max_{y\in\mathcal Y}
		\bigl\{f(x,y)-\langle\lambda,c(x,y)\rangle\bigr\}.
	\end{equation}
	Thus, when the inner feasible set is nonempty for every $x$,
	problem~\eqref{prob:P} has the same optimal value as
	\[
	\inf_{(x,\lambda)\in\mathcal X\times\mathbb R_+^p}
	\max_{y\in\mathcal Y}
	\bigl\{f(x,y)-\langle\lambda,c(x,y)\rangle\bigr\}.
	\]
	This value equivalence does not imply attainment in $\lambda$.
	At a local minimax point of the lifted problem, with $(x,\lambda)$
	as the minimizing block, the first-order necessary conditions are
	\begin{align*}
		0&\in\nabla_x f(x,y)-J_xc(x,y)^\top\lambda+N_{\mathcal X}(x),\\
		0&\in-\nabla_y f(x,y)+J_yc(x,y)^\top\lambda+N_{\mathcal Y}(y),\\
		0&\in-c(x,y)+N_{\mathbb R_+^p}(\lambda).
	\end{align*}
	The last inclusion is equivalent to feasibility, dual nonnegativity,
	and complementarity. These necessary conditions motivate the
	following approximate KKT criterion, which does not assert
	local minimax optimality.

%	{\color{blue}
%		Following~\cite[Section~3.1]{lu23}, we relate the KKT criterion to problem~\eqref{prob:P}. For any fixed $x\in\mathcal X$ with a nonempty inner feasible set, the compactness of $\mathcal Y$, concavity of $f(x,\cdot)$, and convexity of $c_i(x,\cdot)$ allow the convex--concave minimax theorem to yield
%		\begin{align}
%			\max_{\substack{y\in\mathcal Y\\ c(x,y)\le0}} f(x,y)
%			=
%			\inf_{\lambda\in\mathbb R_+^p}\max_{y\in\mathcal Y}
%			\bigl\{f(x,y)-\langle\lambda,c(x,y)\rangle\bigr\}.
%			\label{pro:kkt}
%		\end{align}
%		Consider the resulting lifted minimax formulation \eqref{pro:kkt}, if $(x,y,\lambda)$ is a local minimax point, then $(x,y,\lambda)$ is a KKT point satisfying the conditions:
%		\begin{align*}
%			0&\in\nabla_x f(x,y)-J_xc(x,y)^\top\lambda+N_{\mathcal X}(x),\\
%			0&\in-\nabla_y f(x,y)+J_yc(x,y)^\top\lambda+N_{\mathcal Y}(y),\\
%			c(x,y)&\le0,\qquad \lambda\ge0,\qquad
%			\lambda_i c_i(x,y)=0,\quad i=1,\ldots,p.
%		\end{align*}
%		The last line can equivalently be written as
%		$
%		0\in-c(x,y)+N_{\mathbb R_+^p}(\lambda),
%		$
%		which is the stationarity condition with respect to the multiplier variable. Motivated by this characterization, the following definition relaxes primal stationarity, feasibility, and complementarity while retaining exact dual nonnegativity.}
	
	\begin{definition}[$\varepsilon$-KKT]\label{def:kkt-gs}
		For $\varepsilon>0$, a triple
		$(x,y,\lambda)\in\mathcal X\times\mathcal Y\times\mathbb R_+^p$
		is an $\varepsilon$-KKT point of \eqref{prob:P} if
		\begin{align*}
			\mathcal R_x(x,y,\lambda)
			&:=\dist\!\left(
			0,\nabla_xf(x,y)-J_xc(x,y)^\top\lambda+N_{\mathcal X}(x)
			\right)\le\varepsilon,\\
			\mathcal R_y(x,y,\lambda)
			&:=\dist\!\left(
			0,-\nabla_yf(x,y)+J_yc(x,y)^\top\lambda+N_{\mathcal Y}(y)
			\right)\le\varepsilon,\\
			\mathcal R_c(x,y)
			&:=\|[c(x,y)]_+\|\le\varepsilon,\\
			\mathcal R_{\rm comp}(x,y,\lambda)
			&:=\sum_{i=1}^p|\lambda_i c_i(x,y)|\le\varepsilon.
		\end{align*}
	\end{definition}
	
	The complementarity criterion also gives
	$|\lambda_i c_i(x,y)|\le\mathcal R_{\rm comp}(x,y,\lambda)\le\varepsilon$
	for every $i$.
	To relate the algorithmic updates to these residuals,
	\eqref{def:qgamma}--\eqref{def:Pgamma} yield
	\begin{align}
		\nabla_x\mathcal L_\gamma(x,y,\lambda)
		&=\nabla_xf(x,y)-J_xc(x,y)^\top
		[\lambda+\gamma c(x,y)]_+,\label{eq:gradqx}\\
		\nabla_y\mathcal L_\gamma(x,y,\lambda)
		&=\nabla_yf(x,y)-J_yc(x,y)^\top
		[\lambda+\gamma c(x,y)]_+.\label{eq:gradqy}
	\end{align}
	Thus, the augmented Lagrangian gradients involve the projected multiplier
	used in the output.
	We also recall the output score in Algorithm~\ref{alg:cpdeg}:
	\begin{equation}\label{eq:St}
		S_t:=
		\|\nabla_x\mathcal L_\gamma(x_{t+1},y_{t+1},\lambda_t)
		+\xi_{t+1}\|^2
		+\|-\nabla_y\mathcal L_\gamma(x_{t+1},y_{t+1},\lambda_t)
		+\rho y_{t+1}+n_{t+1}\|^2.
	\end{equation}
		The analysis below bounds this score and the feasibility and complementarity
	residuals to establish the $\varepsilon$-KKT guarantees.
	
	\subsection{Auxiliary Results and KKT Guarantee}\label{sec:aux-kkt}
	We first establish nonnegativity and uniform boundedness of the relaxed multiplier sequence.
	
	\begin{lemma}\label{lem:lambda-state}
		Suppose Assumption~\ref{ass:basic} holds, $\gamma>0$, $\lambda_0=0$,
		$0<\theta_t\le\gamma$, and
		\begin{equation}\label{eq:theta-sum-preliminary}
			\sum_{t=0}^{\infty}\theta_t\le1.
		\end{equation}
		Then Algorithm~\ref{alg:cpdeg} satisfies, for every $t\ge0$,
		\begin{equation}\label{eq:lambda-invariant}
			\lambda_t\in\mathbb R_+^p,\qquad
			\|\lambda_{t+1}-\lambda_t\|\le\theta_tC_c,\qquad
			\|\lambda_t\|\le C_c.
		\end{equation}
	\end{lemma}
	
	\begin{proof}{}
		The relaxed multiplier update is a convex combination of
		$\lambda_t$ and $\widetilde\lambda_{t+1}\in\mathbb R_+^p$.
		Since $\lambda_0=0$, induction gives $\lambda_t\in\mathbb R_+^p$.
		For $a\ge0$ and $s\in\mathbb R$, nonexpansiveness of projection yields
		\begin{equation}\label{eq:scalar-projection-displacement}
			\bigl|[a+\gamma s]_+-a\bigr|\le\gamma|s|.
		\end{equation}
		Applying this inequality componentwise and using (A3), we obtain
		\[
		\|\widetilde\lambda_{t+1}-\lambda_t\|
		\le\gamma\|c(x_{t+1},y_{t+1})\|
		\le\gamma C_c.
		\]
		Hence
		\[
		\|\lambda_{t+1}-\lambda_t\|
		=\frac{\theta_t}{\gamma}
		\|\widetilde\lambda_{t+1}-\lambda_t\|
		\le\theta_tC_c.
		\]
		Summing the increments from $\lambda_0=0$ and using
		\eqref{eq:theta-sum-preliminary} gives
		$\|\lambda_t\|\le C_c\sum_{s<t}\theta_s\le C_c$.
	\end{proof}
	
The next lemma establishes the concavity and Lipschitz bounds of the augmented Lagrangian used in the analysis.

\begin{lemma}\label{lem:penalty}
	Suppose Assumption~\ref{ass:basic} holds, $\gamma\ge1$, and
	$\lambda,\lambda'\in\mathbb R_+^p$ satisfy
	$\max\{\|\lambda\|,\|\lambda'\|\}\le C_c$.
	For every $x\in\mathcal X$, $P_\gamma(x,\cdot,\lambda)$ is convex and
	continuously differentiable, and $\mathcal L_\gamma(x,\cdot,\lambda)$
	is concave. Under Assumption~\ref{ass:strong-y}, the latter is
	$\sigma_y$-strongly concave.
	Moreover, for $u=(x,y)\in\mathcal X\times\mathcal Y$,
	\begin{align}
		|\mathcal L_\gamma(u,\lambda')-\mathcal L_\gamma(u,\lambda)|
		&\le C_c\|\lambda'-\lambda\|,\label{eq:q-lambda-value}\\
		\|\nabla_u\mathcal L_\gamma(u,\lambda')
		-\nabla_u\mathcal L_\gamma(u,\lambda)\|
		&\le J_c\|\lambda'-\lambda\|,\label{eq:q-lambda-gradient}
	\end{align}
	and $\nabla_u\mathcal L_\gamma(u,\lambda)$ is Lipschitz in $u$ on
	$\mathcal X\times\mathcal Y$, uniformly over such $\lambda$, with constant
	\begin{equation}\label{eq:Lgamma}
		L_\gamma^{\rm aug}:=
		L_f+L_{Jc}(C_c+\gamma C_c)+\gamma J_cL_c.
	\end{equation}
\end{lemma}

\begin{proof}{}
	For $a\ge0$, the function
	$\varphi_a(s):=([a+\gamma s]_+^2-a^2)/(2\gamma)$
	is continuously differentiable, convex, and nondecreasing, with
	$\varphi_a'(s)=[a+\gamma s]_+$.
	Since each $c_i(x,\cdot)$ is convex,
	$P_\gamma(x,\cdot,\lambda)
	=\sum_{i=1}^p\varphi_{\lambda_i}(c_i(x,\cdot))$
	is convex and continuously differentiable.
	Subtracting it from $f(x,\cdot)$ preserves concavity and, under
	Assumption~\ref{ass:strong-y}, $\sigma_y$-strong concavity.
	Differentiating \eqref{def:Pgamma} gives
	\[
	\nabla_\lambda P_\gamma(u,\lambda)
	=\frac{[\lambda+\gamma c(u)]_+-\lambda}{\gamma}.
	\]
	By \eqref{eq:scalar-projection-displacement},
	$\|\nabla_\lambda P_\gamma(u,\lambda)\|\le\|c(u)\|\le C_c$.
	Integrating along the segment from $\lambda$ to $\lambda'$ proves
	\eqref{eq:q-lambda-value}.
	Also, $\nabla_uP_\gamma(u,\lambda)
	=Jc(u)^\top[\lambda+\gamma c(u)]_+$, so projection nonexpansiveness gives
	\[
	\|\nabla_uP_\gamma(u,\lambda')
	-\nabla_uP_\gamma(u,\lambda)\|
	\le J_c\|\lambda'-\lambda\|,
	\]
	proving \eqref{eq:q-lambda-gradient}.
	
	Finally, for $u,u'\in\mathcal X\times\mathcal Y$, write
	$w(u):=[\lambda+\gamma c(u)]_+$.
	Then $\|w(u)\|\le C_c+\gamma C_c$ and
	$\|w(u)-w(u')\|\le\gamma L_c\|u-u'\|$.
	Consequently,
	\begin{align*}
		\|\nabla_uP_\gamma(u,\lambda)
		-\nabla_uP_\gamma(u',\lambda)\|
		&\le \|Jc(u)-Jc(u')\|\,\|w(u)\|
		+\|Jc(u')\|\,\|w(u)-w(u')\|\\
		&\le
		\bigl[L_{Jc}(C_c+\gamma C_c)+\gamma J_cL_c\bigr]
		\|u-u'\|.
	\end{align*}
	Combining this bound with the $L_f$-Lipschitz continuity of $\nabla f$
	yields \eqref{eq:Lgamma}.
\end{proof}
Use $\sigma_y>0$ from Assumption~\ref{ass:strong-y} in the strongly concave case
and set $\sigma_y:=0$ in the merely concave case. Define
\begin{equation}\label{eq:ellgamma}
	\mu_y:=\sigma_y+\rho,\qquad
	\ell_\gamma:=\max\{1,L_\gamma^{\mathrm{aug}},\mu_y\}.
\end{equation}
For $0\le\rho\le1$ and $\gamma\ge1$, we have
$\ell_\gamma=\mathcal O(\gamma)$.

The next lemma identifies the projection residuals as normal-cone vectors,
which will be used to bound the stationarity residuals.

\begin{lemma}\label{lem:normals}
	Suppose Assumption~\ref{ass:basic} holds and $h>0$, $\vartheta\ge0$.
	Algorithm~\ref{alg:cpdeg} satisfies
	\[
	\xi_t\in N_{\mathcal X}(x_t),\qquad
	n_t\in N_{\mathcal Y}(y_t),\qquad t=0,\ldots,T.
	\]
\end{lemma}

\begin{proof}{}
	The claim holds at $t=0$ because $\xi_0=n_0=0$.
	By the characterization of Euclidean projection, the correction steps give
	\begin{align*}
		x_t-x_{t+1}
		-h\nabla_xG
		(\widetilde x_t,\widetilde y_t,z_t,\lambda_t)
		&=h\xi_{t+1}\in N_{\mathcal X}(x_{t+1}),\\
		y_t-y_{t+1}
		+h\bigl(\nabla_yG
		(\widetilde x_t,\widetilde y_t,z_t,\lambda_t)+\bar v_t\bigr)
		&=h(1+\vartheta)n_{t+1}\in N_{\mathcal Y}(y_{t+1}).
	\end{align*}
	Since $h>0$, $\vartheta\ge0$, and normal cones are cones,
	the claim follows.
\end{proof}

Condition (A4) also ensures that the coupled constraints admit a feasible
$y\in\mathcal Y$ for every $x\in\mathcal X$.

\begin{lemma}\label{lem:feasible-start}
	Under Assumption~\ref{ass:basic}, for every $x\in\mathcal X$
	there exists $y\in\mathcal Y$ such that $c(x,y)\le0$.
\end{lemma}

\begin{proof}{}
	Fix $x\in\mathcal X$.
	The convex, continuously differentiable function
	$g_x(y):=\frac12\|[c(x,y)]_+\|^2$
	attains its minimum over the compact set $\mathcal Y$ at some $\bar y$.
	First-order optimality gives
	\[
	-J_yc(x,\bar y)^\top[c(x,\bar y)]_+
	\in N_{\mathcal Y}(\bar y).
	\]
	If $g_x(\bar y)>0$, let $d=d(x,\bar y)$ be the direction in (A4).
	Since $d\in N_{\mathcal Y}^*(\bar y)$,
	\[
	\begin{aligned}
		0
		&\ge
		\left\langle d,-J_yc(x,\bar y)^\top[c(x,\bar y)]_+\right\rangle\\
		&=-\sum_{i:c_i(x,\bar y)>0}
		[c_i(x,\bar y)]_+\nabla_yc_i(x,\bar y)^\top d \ge\delta\sum_i[c_i(x,\bar y)]_+>0,
	\end{aligned}
	\]
	a contradiction.
	Thus $g_x(\bar y)=0$, so $c(x,\bar y)\le0$.
\end{proof}		
We next define the regularized saddle value and the potential used in the
remaining analysis.
	Throughout the analysis below, Assumption~\ref{ass:basic} holds, and
	the parameters satisfy $\gamma\ge1$, $0\le\rho\le1$, and
	$0<\mu_y\le\ell_\gamma$.
	Unless otherwise specified, any fixed multiplier
	$\lambda\in\mathbb R_+^p$ satisfies $\|\lambda\|\le C_c$.
	In the strongly concave case, we additionally impose
	Assumption~\ref{ass:strong-y}.
Set
\begin{equation}\label{eq:aux-parameters}
	r_\gamma:=2\ell_\gamma,\qquad
	s_\gamma:=\sqrt{2\ell_\gamma\mu_y}.
\end{equation}
Since $\mathcal L_\gamma(\cdot,y,\lambda)$ is
$\ell_\gamma$-weakly convex, this choice makes
$G(\cdot,y,z,\lambda)$ $\ell_\gamma$-strongly convex.

Fix $\gamma\ge1$, $0\le\rho\le1$ with $\mu_y>0$,
$\lambda\in\mathbb R_+^p$ with $\|\lambda\|\le C_c$,
and $z\in\mathbb R^n$. Define
\[
p_\rho(z,\lambda):=
\min_{x\in\mathcal X}\max_{y\in\mathcal Y}
G(x,y,z,\lambda).
\]
The strong convexity--strong concavity of $G$ ensures a unique
saddle point $(x^\star(z,\lambda),y^\star(z,\lambda))$.
Moreover, $p_\rho(\cdot,\lambda)$ is differentiable, with
\begin{equation}\label{eq:gradp}
	\nabla_z p_\rho(z,\lambda)
	=2\ell_\gamma\bigl(z-x^\star(z,\lambda)\bigr).
\end{equation}

For fixed $(z,\lambda)$, we use an auxiliary energy to control the
stationarity and saddle-tracking errors of a single prediction--correction
step. For $w=(x,y,\xi,n,v)$ with
$\xi\in N_{\mathcal X}(x)$ and $n\in N_{\mathcal Y}(y)$, define
\begin{align}
	\mathcal E_z^\lambda(w)
	&:=-\frac{16\ell_\gamma-s_\gamma}{16\ell_\gamma}
	\bigl(G(x,y,z,\lambda)-p_\rho(z,\lambda)\bigr)
	\notag\\
	&\quad+\frac{1}{\ell_\gamma}
	\|\nabla_xG(x,y,z,\lambda)+\xi\|^2
	+\frac{1}{\ell_\gamma}
	\|-\nabla_yG(x,y,z,\lambda)+n-v\|^2
	\notag\\
	&\quad+\frac{s_\gamma}{16\ell_\gamma}
	\langle v,y-y^\star(z,\lambda)\rangle
	+\frac{\mu_y}{256}\|y-y^\star(z,\lambda)\|^2.
	\label{eq:E}
\end{align}
Its coefficients are chosen to support both the lower bound in
Lemma~\ref{lem:energy-lower} and the one-step dissipation estimate.

The full potential combines the value-function gap, which measures
progress of the center update, with the auxiliary energy, which accounts
for the inexact saddle update:
\begin{equation}\label{eq:V}
	V_t:=
	p_\rho(z_t,\lambda_t)-\inf_z p_\rho(z,\lambda_t)
	+\mathcal E_{z_t}^{\lambda_t}(w_t)
	+\frac{s_\gamma}{256\ell_\gamma^2}\|v_t\|^2.
\end{equation}
The final term controls the sign-indefinite coupling in
$\mathcal E_{z_t}^{\lambda_t}(w_t)$ through square completion and
the contraction of the damping update.

		\begin{lemma}\label{lem:initial-bound}
			Suppose Assumption~\ref{ass:basic} holds and~\eqref{eq:aux-parameters} is satisfied.
			Let
			\[
			f_{\max}:=\max_{\mathcal X\times\mathcal Y}f(x,y),
			\qquad
			f_{\min}:=\min_{\mathcal X\times\mathcal Y}f(x,y),
			\]
			and suppose $\mu_y>0$. Initialize Algorithm~\ref{alg:cpdeg} with
			$c(x_0,y_0)\le0$, $\lambda_0=0$, $z_0=x_0$, and $\xi_0=n_0=v_0=0$.
			If $0\le\rho\le1$, then
			\begin{equation}\label{eq:Bsig}
				V_0\le
				B_{\sigma_y}:=
				2(f_{\max}-f_{\min})+G_x^2+2G_y^2
				+\left(3+\frac{\sigma_y+1}{64}\right)D_Y^2.
			\end{equation}
			%{\color{blue}In particular, $B_{\sigma_y}$ is independent of $\gamma$ and $\rho$.}
		\end{lemma}
		
		\begin{proof}{Proof.}
			\noindent
			Define $\phi_{0,\rho}(x):=
			\max_{y\in\mathcal Y}
			\left(\mathcal L_\gamma(x,y,0)-\frac{\rho}{2}\|y\|^2\right)$.
			Lemma~\ref{lem:feasible-start} implies
			$\phi_{0,\rho}(x)\ge f_{\min}-\rho D_Y^2/2$ for every $x$, whereas
			$\mathcal L_\gamma(\cdot,\cdot,0)\le f$ gives $\phi_{0,\rho}(x)\le f_{\max}$. Since
			\[
			\inf_z p_\rho(z,0)=\inf_{x\in\mathcal X}\phi_{0,\rho}(x),
			\]
			we obtain
			\[
			p_\rho(x_0,0)-\inf_z p_\rho(z,0)
			\le f_{\max}-f_{\min}+\frac{\rho D_Y^2}{2}.
			\]
			At the feasible initial point,
			$P_\gamma(x_0,y_0,0)=0$ and $\nabla_u P_\gamma(x_0,y_0,0)=0$. Therefore
			\[
			p_\rho(x_0,0)-G(x_0,y_0,x_0,0)
			\le f_{\max}-f_{\min}+\frac{\rho D_Y^2}{2},
			\]
			\[
			\nabla_xG(x_0,y_0,x_0,0)=\nabla_xf(x_0,y_0),\qquad
			-\nabla_yG(x_0,y_0,x_0,0)
			=-\nabla_yf(x_0,y_0)+\rho y_0.
			\]
			Substituting the last displays into
			\eqref{eq:E}--\eqref{eq:V}, using $\ell_\gamma\ge1$,
			$\|y_0-y^\star(x_0,0)\|\le 2D_Y$, and
			$\|a+b\|^2\le2\|a\|^2+2\|b\|^2$, gives
			\[
			V_0
			\le2(f_{\max}-f_{\min})+G_x^2+2G_y^2
			+\left(\rho+2\rho^2+\frac{\sigma_y+\rho}{64}\right)D_Y^2.
			\]
			Finally $0\le\rho\le1$ proves \eqref{eq:Bsig}.
		\end{proof}
		
We next derive two estimates that convert augmented stationarity bounds into KKT guarantees. The first bounds the effective multiplier under the extended MFCQ condition.
		\begin{lemma}\label{lem:mu-bound}
			Suppose Assumption~\ref{ass:basic} holds and $\gamma>0$, $\rho\ge0$, $r\ge0$, $(x,y)\in\mathcal X\times\mathcal Y$. Let
			$\lambda\in\mathbb R_+^p$ satisfy $\|\lambda\|\le C_c$, set
			\begin{equation}\label{eq:effective-mu}
				\mu=[\lambda+\gamma c(x,y)]_+,
			\end{equation}
			and suppose
$\dist\!\left(0,-\nabla_yf(x,y)+J_yc(x,y)^\top\mu
			+\rho y+N_{\mathcal Y}(y)\right)\le r.$
			Then
			\begin{equation}\label{eq:mu-bound}
				\|\mu\|_1\le \sqrt p\,C_c
				+\frac{Z}{\delta}\bigl(r+G_y+\rho D_Y+J_cC_c\bigr).
			\end{equation}
			If $r\le1$ and $0\le\rho\le1$, then
			\begin{equation}\label{eq:Mmu}
				\|\mu\|_1\le M_{\rm aug}:=
				\sqrt p\,C_c+
				\frac{Z}{\delta}(1+G_y+D_Y+J_cC_c).
			\end{equation}
		\end{lemma}
		
		\begin{proof}{Proof.}
			\noindent
			There are $n\in N_{\mathcal Y}(y)$ and $e$ with $\|e\|\le r$ such that
			\[
			-\nabla_yf(x,y)+J_yc(x,y)^\top\mu+\rho y+n=e.
			\]
			Let $d=d(x,y)$ be supplied by Assumption~\ref{ass:basic}(A4). Since
			$d\in N_{\mathcal Y}^*(y)$, $\langle n,d\rangle\le0$, and hence
			\[
			\langle J_yc(x,y)^\top\mu,d\rangle
			=\langle e+\nabla_yf(x,y)-\rho y-n,d\rangle
			\ge-Z(r+G_y+\rho D_Y).
			\]
			Let $I_+:=\{i:c_i(x,y)>0\}$ and $I_0:=\{1,\ldots,p\}\setminus I_+$. For $i\in I_0$,
			$0\le\mu_i=[\lambda_i+\gamma c_i(x,y)]_+\le\lambda_i$. Hence
			\begin{equation}\label{eq:mu-nonviolated}
				\|\mu_{I_0}\|\le C_c,
				\qquad
				\|\mu_{I_0}\|_1\le\sqrt p\,C_c.
			\end{equation}
			On the violated indices, Assumption~\ref{ass:basic}(A4) gives
			\[
			\sum_{i\in I_+}\mu_i\nabla_yc_i(x,y)^\top d
			\le-\delta\|\mu_{I_+}\|_1.
			\]
			The remaining contribution satisfies
			\[
			\left|\sum_{i\in I_0}\mu_i\nabla_yc_i(x,y)^\top d\right|
			\le ZJ_c\|\mu_{I_0}\|\le ZJ_cC_c.
			\]
			Combining the last three displays yields
			\[
			\|\mu_{I_+}\|_1
			\le\frac Z\delta(r+G_y+\rho D_Y+J_cC_c).
			\]
			Adding \eqref{eq:mu-nonviolated} proves \eqref{eq:mu-bound},
			\eqref{eq:Mmu} follows from $r\le1$ and $\rho\le1$.
		\end{proof}
		
	The next lemma translates the effective multiplier bound into feasibility and complementarity bounds.
		\begin{lemma}\label{lem:feas-comp}
			Let $\gamma>0$, $(x,y)\in\mathcal X\times\mathcal Y$, $\lambda\in\mathbb R_+^p$, $\mu$ be given by
			\eqref{eq:effective-mu}, and suppose
			$\|\lambda\|\le C_c$ and $\|\mu\|_1\le M$. Then
			\begin{equation}\label{eq:feas-comp}
				\|[c(x,y)]_+\|\le\frac{M}{\gamma},
				\qquad
				\sum_{i=1}^p|\mu_i c_i(x,y)|
				\le\frac{M(M+C_c)}{\gamma}.
			\end{equation}
		\end{lemma}
		
		\begin{proof}{Proof.}
			\noindent
			If $c_i(x,y)>0$, then
			$\mu_i=\lambda_i+\gamma c_i(x,y)\ge\gamma c_i(x,y)$; hence
			\[
			\|[c(x,y)]_+\|\le\frac{\|\mu\|}{\gamma}
			\le\frac M\gamma.
			\]
			If $\mu_i>0$, the projection identity gives
			$\gamma c_i(x,y)=\mu_i-\lambda_i$. Consequently, for every $i$,
			\[
			\mu_i c_i(x,y)=\frac{\mu_i(\mu_i-\lambda_i)}{\gamma}.
			\]
			Taking absolute values before summation, using $\mu,\lambda\ge0$, and applying Cauchy-Schwarz gives
			\begin{align*}
				\sum_{i=1}^p|\mu_i c_i(x,y)|
				&=\frac1\gamma\sum_{i=1}^p\mu_i|\mu_i-\lambda_i|\le\frac{\|\mu\|^2+\mu^\top\lambda}{\gamma}
				\le\frac{\|\mu\|(\|\mu\|+\|\lambda\|)}{\gamma}\le\frac{M(M+C_c)}{\gamma},
			\end{align*}
			where the last step uses $\|\lambda\|\le C_c$ and
			$\|\mu\|\le\|\mu\|_1\le M$.
		\end{proof}
		
		\begin{lemma}\label{lem:energy-lower}
			Suppose Assumption~\ref{ass:basic} holds,~\eqref{eq:aux-parameters} is satisfied and $\xi\in N_{\mathcal X}(x)$, $n\in N_{\mathcal Y}(y)$.
			For a fixed multiplier and center, suppress the last two arguments of $G$,
			set $p=p_\rho(z,\lambda)$, and write
			$x^\star:=x^\star(z,\lambda)$ and $y^\star:=y^\star(z,\lambda)$;
			$a=\nabla_xG(x,y)+\xi$, and $b=-\nabla_yG(x,y)+n$.
			Then
			\begin{align}\label{eq:energy-lower}
				\mathcal E_z^\lambda(w)
				&\ge\frac{\|a\|^2}{2\ell_\gamma}
				+\frac{s_\gamma}{32}\|x-x^\star\|^2
				+\frac{\mu_y}{2}\|y-y^\star\|^2+\frac1{\ell_\gamma}
				\left\|b-v-\frac{s_\gamma}{32}(y-y^\star)\right\|^2\ge0.
			\end{align}
		\end{lemma}
		
		\begin{proof}{Proof.}
			Strong convexity and $\xi\in N_{\mathcal X}(x)$ give
			\[
			G(\bar x,y)\ge G(x,y)+\langle a,\bar x-x\rangle+\frac{\ell_\gamma}{2}\|\bar x-x\|^2
			\ge G(x,y)-\frac{\|a\|^2}{2\ell_\gamma}.
			\]
			Taking $\bar x=x^\star$ and then using strong concavity of
			$G(x^\star,\cdot)$ at $y^\star$ gives
			\[
			G(x,y)-p\le\frac{\|a\|^2}{2\ell_\gamma}
			-\frac{\mu_y}{2}\|y-y^\star\|^2.
			\]
			Likewise, $n\in N_{\mathcal Y}(y)$ and strong concavity imply
			\[
			\langle b,y-y^\star\rangle
			\ge G(x,y^\star)-G(x,y)+\frac{\mu_y}{2}\|y-y^\star\|^2,
			\qquad
			G(x,y^\star)-p\ge\frac{\ell_\gamma}{2}\|x-x^\star\|^2.
			\]
			Since $s_\gamma^2=2\ell_\gamma\mu_y$, we get
			\begin{align*}
				\mathcal E_z^\lambda(w)
				&=-\left(1-\frac{s_\gamma}{16\ell_\gamma}\right)(G(x,y)-p)
				+\frac{\|a\|^2}{\ell_\gamma}
				+\frac{s_\gamma}{16\ell_\gamma}\langle b,y-y^\star\rangle\\
				&\quad+\frac1{\ell_\gamma}
				\left\|b-v-\frac{s_\gamma}{32}(y-y^\star)\right\|^2
				+\frac{\mu_y}{512}\|y-y^\star\|^2.
			\end{align*}
			Substituting the preceding inequalities gives
			\begin{align*}
				\mathcal E_z^\lambda(w)
				&\ge -(G(x,y)-p)+\frac{\|a\|^2}{\ell_\gamma}
				+\frac{s_\gamma}{16\ell_\gamma}\bigl(G(x,y^\star)-p\bigr)\\
				&\quad+\left(\frac{s_\gamma\mu_y}{32\ell_\gamma}
				+\frac{\mu_y}{512}\right)\|y-y^\star\|^2
				+\frac1{\ell_\gamma}
				\left\|b-v-\frac{s_\gamma}{32}(y-y^\star)\right\|^2\\
				&\ge\frac{\|a\|^2}{2\ell_\gamma}
				+\frac{s_\gamma}{32}\|x-x^\star\|^2
				+\frac{\mu_y}{2}\|y-y^\star\|^2
				+\frac1{\ell_\gamma}
				\left\|b-v-\frac{s_\gamma}{32}(y-y^\star)\right\|^2,
			\end{align*}
			which proves \eqref{eq:energy-lower}.
		\end{proof}
		
	The next lemma establishes potential decrease for a fixed multiplier and bounds the error induced by the center update.
		\begin{lemma}\label{lem:technical}
			Suppose Assumption~\ref{ass:basic} holds and~\eqref{eq:aux-parameters} is satisfied. Consider an iteration of Algorithm~\ref{alg:cpdeg} with $\lambda_t\ge0$, $\|\lambda_t\|\le C_c$.
			Set
			\begin{equation}\label{eq:solver-parameters}
				L_G:=3\ell_\gamma+\rho,\ 
				h:=\frac{1}{64L_G},\ 
				\alpha:=\left(1+\frac{hs_\gamma}{16}\right)^{-1},\  \vartheta:=\frac{\alpha h(8\ell_\gamma-s_\gamma)}{16},  \beta:=\frac{hs_\gamma}{4096}.
			\end{equation}
%			\begin{equation}\label{eq:solver-parameters-2}
%				\vartheta:=\frac{\alpha h(8\ell_\gamma-s_\gamma)}{16},
%				\qquad \beta:=\frac{hs_\gamma}{4096}.
%			\end{equation}
			Then $s_\gamma\le\sqrt2\ell_\gamma$, $0<\alpha<1$,
			$\vartheta\ge0$, and $0<\beta<1$.
			For iteration $t$, write $G_t(x,y):=G(x,y,z_t,\lambda_t)$, $(x_t^\star,y_t^\star):=
			(x^\star(z_t,\lambda_t),y^\star(z_t,\lambda_t)).$
			Define $A_{x,t} :=\nabla_xG_t(x_{t+1},y_{t+1})+\xi_{t+1}$, $A_{y,t}:=-\nabla_yG_t(x_{t+1},y_{t+1})+n_{t+1}-v_{t+1}$,
			\begin{align}
				I_t
				&:=\|\nabla_xG_t(x_{t+1},y_{t+1})+\xi_{t+1}
				-\nabla_xG_t(x_t,y_t)-\xi_t\|^2\notag\\
				&\quad+\|-\nabla_yG_t(x_{t+1},y_{t+1})+n_{t+1}-v_{t+1}
				+\nabla_yG_t(x_t,y_t)-n_t+v_t\|^2,\label{eq:It}\\
%				A_{x,t}
%				&:=\nabla_xG_t(x_{t+1},y_{t+1})+\xi_{t+1},
%				\qquad
%				A_{y,t}:=-\nabla_yG_t(x_{t+1},y_{t+1})
%				+n_{t+1}-v_{t+1},\notag\\
				\mathcal D_t
				&:=\frac h4\|A_{x,t}\|^2
				+\frac{3hs_\gamma}{512\ell_\gamma}\|A_{y,t}\|^2
				+\frac{7hs_\gamma}{4096\ell_\gamma}\|v_{t+1}\|^2\notag\\
				&\quad+\frac{h\mu_y s_\gamma}{128}
				\|y_{t+1}-y_t^\star\|^2
				+\frac{hs_\gamma(16\ell_\gamma-s_\gamma)}{1024}
				\|x_{t+1}-x_t^\star\|^2
				+\frac{I_t}{4\ell_\gamma},\label{eq:Dt}\\
				\mathcal B_t^\lambda
				&:=\frac{hs_\gamma}{8\ell_\gamma}
				\|\nabla_xG(x_{t+1},y_{t+1},z_t,\lambda)
				+\xi_{t+1}\|^2
				+\frac{7hs_\gamma}{4096\ell_\gamma}\|v_{t+1}\|^2\notag\\
				&\quad+\frac{h\mu_y s_\gamma}{128}
				\|y_{t+1}-y^\star(z_t,\lambda)\|^2
				+\frac{7h\ell_\gamma s_\gamma}{512}
				\|x_{t+1}-x^\star(z_t,\lambda)\|^2.\label{eq:Bt}
			\end{align}
			These quantities satisfy
			\begin{equation}\label{eq:B-le-D}
				0\le\mathcal B_t^{\lambda_t}\le\mathcal D_t.
			\end{equation}
			At fixed $(z_t,\lambda_t)$, the energy satisfies
			\begin{align}
				&\mathcal E_{z_t}^{\lambda_t}(w_{t+1})
				-\mathcal E_{z_t}^{\lambda_t}(w_t)
				+\frac{s_\gamma}{256\ell_\gamma^2}
				\bigl(\|v_{t+1}\|^2-\|v_t\|^2\bigr)\notag\\
				&\quad+\frac{hs_\gamma}{32}
				\left(\mathcal E_{z_t}^{\lambda_t}(w_{t+1})
				+\frac{s_\gamma}{256\ell_\gamma^2}\|v_{t+1}\|^2\right)
				\le-\mathcal D_t.\label{eq:fixed-center}
			\end{align}
			For any fixed $\lambda\in\mathbb R_+^p$ with $\|\lambda\|\le C_c$,
			the center update satisfies
			\begin{align}
				\mathcal E_{z_{t+1}}^\lambda(w_{t+1})
				-\mathcal E_{z_t}^\lambda(w_{t+1})
				&\le\frac14\mathcal B_t^\lambda
				+\frac{352\beta^2}{\ell_\gamma hs_\gamma}
				\|\nabla_z p_\rho(z_t,\lambda)\|^2\notag\\
				&\quad+\frac{1408\ell_\gamma\beta^2}{hs_\gamma}
				\|x_{t+1}-x^\star(z_t,\lambda)\|^2.
				\label{eq:center-motion}
			\end{align}
		\end{lemma}

		\begin{proof}{Proof.}
			Throughout iteration $t$, all gradients are taken with respect to $G_t$.
			Define
			\begin{align*}
				e_{x,t}
				&:=\nabla_xG_t(\widetilde x_t,\widetilde y_t)
				-\nabla_xG_t(x_{t+1},y_{t+1}),\\
				e_{y,t}
				&:=(1+\vartheta)
				\bigl(\nabla_yG_t(x_{t+1},y_{t+1})
				-\nabla_yG_t(\widetilde x_t,\widetilde y_t)\bigr).
			\end{align*}
			The correction updates give
			\begin{align}
				x_{t+1}-x_t
				&=-h\bigl(\nabla_xG_t(x_{t+1},y_{t+1})
				+\xi_{t+1}+e_{x,t}\bigr),\label{eq:error-x}\\
				y_{t+1}-y_t
				&=-h\bigl(-\nabla_yG_t(x_{t+1},y_{t+1})
				+n_{t+1}-v_{t+1}+e_{y,t}\bigr).\label{eq:error-y}
			\end{align}
			Projection nonexpansiveness and $L_G$-Lipschitz continuity give
			\begin{align*}
			\|(\widetilde x_t-x_{t+1},\widetilde y_t-y_{t+1})\|
			&\le h\left(\sqrt{I_t}+\sqrt{\|e_{x,t}\|^2+\|e_{y,t}\|^2}\right),\\
			\sqrt{\|e_{x,t}\|^2+\|e_{y,t}\|^2}
			&\le(1+\vartheta)L_G\|(\widetilde x_t-x_{t+1},\widetilde y_t-y_{t+1})\|.
			\end{align*}
			Since $hL_G=1/64$ and $0\le\vartheta\le h\ell_\gamma/2\le1/128$, absorbing the error term yields
			\begin{align}\label{eq:error-bound}
				\sqrt{\|e_{x,t}\|^2+\|e_{y,t}\|^2}
				\le\frac{h(1+\vartheta)L_G}{1-h(1+\vartheta)L_G}\sqrt{I_t}
				\le2hL_G\sqrt{I_t}.
			\end{align}
			Since $s_\gamma\le2\ell_\gamma$,
			$hs_\gamma/(8\ell_\gamma)\le h/4$ and
			$7h\ell_\gamma s_\gamma/512
			\le hs_\gamma(16\ell_\gamma-s_\gamma)/1024$.
			Comparing \eqref{eq:Bt} and \eqref{eq:Dt} proves
			\eqref{eq:B-le-D}.
			
			We first prove the fixed-center estimate. Write
			$p_i:=\nabla_xG_t(x_i,y_i)+\xi_i$ and
			$b_i:=-\nabla_yG_t(x_i,y_i)+n_i$ for $i=t,t+1$,
			and set $P:=p_{t+1}$, $R:=b_{t+1}-v_{t+1}$,
			$X:=x_{t+1}-x_t^\star$, $Y:=y_{t+1}-y_t^\star$,
			$\Delta x:=x_{t+1}-x_t$, $\Delta y:=y_{t+1}-y_t$,
			$\Delta v:=v_{t+1}-v_t$, $U:=p_{t+1}-p_t$ and $W:=b_{t+1}-v_{t+1}-b_t+v_t$.
			Thus $I_t=\|U\|^2+\|W\|^2$, $P=A_{x,t}$, and $R=A_{y,t}$.
			The updates and the definitions of $\alpha$ and $\vartheta$ imply
			\begin{equation}\label{eq:frozen-dynamics}
				\Delta x=-h(P+e_{x,t}),\qquad
				\Delta y=-h(R+e_{y,t}),\qquad
				\Delta v=-\frac{h\ell_\gamma}{2}b_{t+1}
				+\frac{hs_\gamma}{16}R.
			\end{equation}
			Strong convexity--strong concavity and normal-cone monotonicity give
			\begin{align}
			&	G_t(x_{t+1},y_{t+1})-G_t(x_t,y_t)
				\ge\langle p_t,\Delta x\rangle
				-\langle b_{t+1},\Delta y\rangle
				+\frac{\ell_\gamma}{2}\|\Delta x\|^2
				+\frac{\mu_y}{2}\|\Delta y\|^2,
				\label{eq:frozen-mixed}\\
				& \langle U,\Delta x\rangle
				+\langle W+\Delta v,\Delta y\rangle
				\ge\ell_\gamma\|\Delta x\|^2
				+\mu_y\|\Delta y\|^2,\label{eq:frozen-monotone}\\
				& \langle b_{t+1},Y\rangle
				+G_t(x_{t+1},y_{t+1})-p_\rho(z_t,\lambda_t)
				\ge\frac{\ell_\gamma}{2}\|X\|^2
				+\frac{\mu_y}{2}\|Y\|^2.\label{eq:frozen-growth}
			\end{align}
			Expanding \eqref{eq:E} with
			$\|a\|^2-\|b\|^2=2\langle a,a-b\rangle-\|a-b\|^2$
			and applying \eqref{eq:frozen-dynamics}--\eqref{eq:frozen-growth}
			and discarding the nonpositive term $-(256\ell_\gamma^2-s_\gamma^2)\|\Delta x\|^2/(2048h)$ yields
			\begin{equation}\label{eq:frozen-decomposition}
				\frac{\ell_\gamma}{2h}
				\left[
				\mathcal E_{z_t}^{\lambda_t}(w_{t+1})
				-\mathcal E_{z_t}^{\lambda_t}(w_t)
				+\frac{hs_\gamma}{32}
				\mathcal E_{z_t}^{\lambda_t}(w_{t+1})
				\right]\le Q_t+H_t,
			\end{equation}
			where
				\begin{align*}
					Q_t&=\frac{32\ell_\gamma-s_\gamma}{64}\|P\|^2-\ell_\gamma\|P+e_{x,t}\|^2-\mu_y\|R+e_{y,t}\|^2-\frac{s_\gamma}{64}\|R\|^2\\
					&\quad+\frac{16\ell_\gamma-s_\gamma}{32}\langle P,e_{x,t}\rangle-\frac{s_\gamma}{32}\langle R,e_{y,t}\rangle
					-\frac{\ell_\gamma\mu_y}{256}\langle Y,e_{y,t}\rangle-\frac{\ell_\gamma\mu_y}{512}\langle R,Y\rangle\\
					&\quad+\left(\frac{\ell_\gamma\mu_y s_\gamma}{16384}-\frac{\mu_y s_\gamma(16\ell_\gamma-s_\gamma)}{2048}\right)\|Y\|^2-\frac{\ell_\gamma s_\gamma(16\ell_\gamma-s_\gamma)}{2048}\|X\|^2,\\
					H_t&=-\frac1{2h}\left\|U-\frac{16\ell_\gamma-s_\gamma}{32}\Delta x\right\|^2-\frac1{2h}\|W\|^2
					-\left(\frac{\mu_y(16\ell_\gamma-s_\gamma)}{64h}+\frac{\ell_\gamma\mu_y}{512h}\right)\|\Delta y\|^2\\
					&\quad-\frac{s_\gamma}{32h}\langle\Delta v,\Delta y\rangle-\frac1h\langle e_{x,t},U\rangle-\frac1h\langle e_{y,t},W\rangle.
			\end{align*}
			To bound $Q_t$, expand the squared $P$-term and use
			$2\ell_\gamma\|P\|\|e_{x,t}\|\le\ell_\gamma\|P\|^2/4+4\ell_\gamma\|e_{x,t}\|^2$.
			The remaining cross terms can be bounded as follows:
			\begin{align*}
			\frac{s_\gamma}{32}|\langle R,e_{y,t}\rangle|&\le\frac{s_\gamma}{256}\|R\|^2+\frac{s_\gamma}{16}\|e_{y,t}\|^2,\quad
			\frac{\ell_\gamma\mu_y}{512}|\langle R,Y\rangle|\le\frac{s_\gamma}{256}\|R\|^2+\frac{\ell_\gamma\mu_y s_\gamma}{8192}\|Y\|^2,\\
			\frac{\ell_\gamma\mu_y}{256}|\langle Y,e_{y,t}\rangle|&\le\frac{\ell_\gamma\mu_y s_\gamma}{16384}\|Y\|^2+\frac{s_\gamma}{32}\|e_{y,t}\|^2.
			\end{align*}
			Using $s_\gamma^2=2\ell_\gamma\mu_y$ and $s_\gamma\le2\ell_\gamma$ gives
			\begin{align}
				Q_t
				&\le-\frac{\ell_\gamma}{4}\|P\|^2
				-\frac{s_\gamma}{128}\|R\|^2
				-\frac{\ell_\gamma\mu_y s_\gamma}{256}\|Y\|^2
				-\frac{\ell_\gamma s_\gamma(16\ell_\gamma-s_\gamma)}
				{2048}\|X\|^2\notag\\
				&\quad+4\ell_\gamma
				\bigl(\|e_{x,t}\|^2+\|e_{y,t}\|^2\bigr).
				\label{eq:current-block}
			\end{align}
			For the increment terms, use
			$\|a-b\|^2\ge\frac12\|a\|^2-\|b\|^2$ and
			\[
			\Delta v
			=-\nabla_yG_t(x_{t+1},y_{t+1})
			+\nabla_yG_t(x_t,y_t)+n_{t+1}-n_t-W.
			\]
			The first two squares in $H_t$ are bounded by
			$
			-\frac{I_t}{4h}+\frac{h(16\ell_\gamma-s_\gamma)^2}{1024}(\|P\|^2+\|e_{x,t}\|^2).
			$
			By normal-cone monotonicity and the $L_G$-Lipschitz continuity of the gradient, the term involving $\Delta v$ is bounded by
			\begin{align*}
				-\frac{s_\gamma}{32h}\langle\Delta v,\Delta y\rangle
				&\le \frac{hL_Gs_\gamma}{32}\bigl(\|P\|^2+\|e_{x,t}\|^2\bigr)\\
				&\quad+\left(\frac{3hL_Gs_\gamma}{32}
				+\frac{h\ell_\gamma\mu_y}{32}\right)
				\bigl(\|R\|^2+\|e_{y,t}\|^2\bigr)
				+\frac{I_t}{32h}.
			\end{align*}
			Moreover, by~\eqref{eq:error-bound}, the last two terms in $H_t$ are bounded above by $2L_G I_t$. Substituting these estimates and using $\mu_y\le s_\gamma$, $\ell_\gamma\le L_G$, and $hL_G=1/64$ yields
			\begin{equation}\label{eq:increment-block}
				H_t\le
				\frac{5\ell_\gamma hL_G}{16}\|P\|^2
				+\frac{hL_Gs_\gamma}{8}\|R\|^2
				+\frac{\ell_\gamma}{2}
				\bigl(\|e_{x,t}\|^2+\|e_{y,t}\|^2\bigr)
				-\frac{7/32-2hL_G}{h}I_t.
			\end{equation}
			Also, \eqref{eq:frozen-dynamics} implies
			$v_t=(1+h\ell_\gamma/2)v_{t+1}
			+h(8\ell_\gamma-s_\gamma)R/16$ and
			$\|v_t\|^2-\|v_{t+1}\|^2
			\ge(h\ell_\gamma/2)
			(\|v_{t+1}\|^2-\|R\|^2)$.
			Hence the additional momentum terms satisfy
			\begin{align*}
				&\frac{s_\gamma}{256\ell_\gamma^2}
				\bigl(\|v_{t+1}\|^2-\|v_t\|^2\bigr)
				+\frac{hs_\gamma}{32}
				\frac{s_\gamma}{256\ell_\gamma^2}\|v_{t+1}\|^2
				\le 
				-\frac{7hs_\gamma}{4096\ell_\gamma}\|v_{t+1}\|^2
				+\frac{hs_\gamma}{512\ell_\gamma}\|R\|^2.
			\end{align*}
			Combining this bound with
			\eqref{eq:frozen-decomposition}--\eqref{eq:increment-block},
			and using \eqref{eq:error-bound} and $hL_G=1/64$,
			proves \eqref{eq:fixed-center}.
			
			We next fix $\lambda$ and estimate the center motion. Write
			$u_z^\star:=(x_z^\star,y_z^\star)
			:=(x^\star(z,\lambda),y^\star(z,\lambda))$ and
			$F_z(u):=\bigl(\nabla_xG(x,y,z,\lambda),
			-\nabla_yG(x,y,z,\lambda)\bigr).$
			Adding the saddle variational inequalities at $z$ and $z'$
			and using strong monotonicity gives
			\begin{align*}
				\ell_\gamma\|x_{z'}^\star-x_z^\star\|^2
				+\mu_y\|y_{z'}^\star-y_z^\star\|^2
				&\le2\ell_\gamma
				\langle z'-z,x_{z'}^\star-x_z^\star\rangle\\
				&\le\frac{\ell_\gamma}{2}
				\|x_{z'}^\star-x_z^\star\|^2
				+2\ell_\gamma\|z'-z\|^2.
			\end{align*}
			Consequently,
			$\|x_{z'}^\star-x_z^\star\|\le2\|z'-z\|$ and
			$\|y_{z'}^\star-y_z^\star\|
			\le\sqrt{2\ell_\gamma/\mu_y}\|z'-z\|$,
			and \eqref{eq:gradp} implies that
			$\nabla_zp_\rho(\cdot,\lambda)$ is $6\ell_\gamma$-Lipschitz.		
			Set $w=w_{t+1}$ and $\zeta=z_{t+1}-z_t$.
			Using the preceding bounds and
			\[
			G(x,y,z+\zeta,\lambda)-G(x,y,z,\lambda)
			=2\ell_\gamma\langle z-x,\zeta\rangle
			+\ell_\gamma\|\zeta\|^2
			\]
			in \eqref{eq:E}, we obtain
			\begin{align*}
				\mathcal E_{z_{t+1}}^\lambda(w)
				-\mathcal E_{z_t}^\lambda(w)
				&\le\Bigl[
				2\ell_\gamma\|x_{t+1}-x^\star(z_t,\lambda)\|+4\|\nabla_xG
				(x_{t+1},y_{t+1},z_t,\lambda)+\xi_{t+1}\|\\
				&\quad
				+\frac18\|v_{t+1}\|+\frac{s_\gamma}{128}
				\|y_{t+1}-y^\star(z_t,\lambda)\|
				\Bigr]\|\zeta\|
				+\frac{1025\ell_\gamma}{128}\|\zeta\|^2\\
				&\le\frac14\mathcal B_t^\lambda
				+\frac{704\ell_\gamma}{hs_\gamma}\|\zeta\|^2.
			\end{align*}
			The last inequality follows from weighted Young inequalities
			and $hs_\gamma\le32$.
			Finally, the center update and \eqref{eq:gradp} give
			\[
			\|\zeta\|^2
			\le2\beta^2\|x_{t+1}-x^\star(z_t,\lambda)\|^2
			+\frac{\beta^2}{2\ell_\gamma^2}
			\|\nabla_zp_\rho(z_t,\lambda)\|^2,
			\]
			proving \eqref{eq:center-motion}.
		\end{proof}
		
		\begin{lemma}\label{lem:lambda-drift}
			Suppose Assumption~\ref{ass:basic} and the conditions of Lemma~\ref{lem:technical} hold.
			Let 
			\begin{align}
				K_E&:=J_c^2\left[
				\frac1{\ell_\gamma}\left(1+\frac{64}{h\ell_\gamma}\right)
				+\frac1{\ell_\gamma}\left(1+\frac{8192}{3hs_\gamma}\right)
				+\frac{128}{7hs_\gamma\mu_y}
				+\frac{1}{32hs_\gamma\mu_y}+\frac{1}{256\mu_y}
				\right],\label{eq:KE}\\
				K_B&:=5J_c^2\left[
				\frac{hs_\gamma}{8\ell_\gamma}
				+\frac{hs_\gamma}{128\mu_y}
				+\frac{7h\ell_\gamma s_\gamma}{512\mu_y^2}
				\right],\label{eq:KB}
			\end{align}
			$\lambda,\lambda'\in\mathbb R_+^p$ satisfy
			$\max\{\|\lambda\|,\|\lambda'\|\}\le C_c$, and put
			$d=\|\lambda'-\lambda\|$. At every fixed center $z$,
			\begin{equation}\label{eq:saddle-lambda-sensitivity}
				\left\|(x^\star(z,\lambda'),y^\star(z,\lambda'))
				-(x^\star(z,\lambda),y^\star(z,\lambda))\right\|
				\le\frac{J_c}{\mu_y}d.
			\end{equation}
			At $z=z_t$, $w=w_{t+1}$, and $\lambda=\lambda_t$,
			\begin{align}
				&\left[p_\rho(z,\lambda')-\inf_{\bar z}p_\rho(\bar z,\lambda')
				+\mathcal E_z^{\lambda'}(w)
				+\frac{s_\gamma}{256\ell_\gamma^2}\|v_{t+1}\|^2\right]\notag\\
				&\quad-
				\left[p_\rho(z,\lambda)-\inf_{\bar z}p_\rho(\bar z,\lambda)
				+\mathcal E_z^\lambda(w)
				+\frac{s_\gamma}{256\ell_\gamma^2}\|v_{t+1}\|^2\right]
				\le\frac1{16}\mathcal D_t+4C_cd+K_Ed^2,
				\label{eq:E-lambda-switch}
			\end{align}
			and
			\begin{equation}\label{eq:B-lambda-switch}
				\mathcal B_t^{\lambda'}
				\le\frac54\mathcal B_t^\lambda+K_Bd^2.
			\end{equation}
		\end{lemma}
		
		\begin{proof}{Proof.}
			Let $u^\star(z,\lambda)=(x^\star(z,\lambda),y^\star(z,\lambda))$.
			For fixed $z$, let $u=(x,y)$ and define
			\[
			F_\lambda(u):=(\nabla_xG(x,y,z,\lambda),-\nabla_yG(x,y,z,\lambda)).
			\]
			Using strong monotonicity gives
			\[
			\mu_y\|u^\star(z,\lambda')-u^\star(z,\lambda)\|^2
			\le\langle F_\lambda(u^\star(z,\lambda'))
			-F_{\lambda'}(u^\star(z,\lambda')),
			u^\star(z,\lambda')-u^\star(z,\lambda)\rangle.
			\]
			\eqref{eq:q-lambda-gradient} and Cauchy--Schwarz prove
			\eqref{eq:saddle-lambda-sensitivity}.
			\eqref{eq:q-lambda-value} also gives
			\[
			|p_\rho(z,\lambda')-p_\rho(z,\lambda)|\le C_cd,
			\qquad
			\left|\inf_z p_\rho(z,\lambda')-\inf_z p_\rho(z,\lambda)\right|
			\le C_cd.
			\]
			At the common point $(x_{t+1},y_{t+1},z_t)$, suppress the first three
			arguments of $G$. The two residuals decompose as
			\begin{align*}
				\nabla_xG(\lambda')+\xi_{t+1}
				&=\nabla_xG(\lambda)+\xi_{t+1}
				+\bigl(\nabla_xG(\lambda')-\nabla_xG(\lambda)\bigr),\notag\\
				-\nabla_yG(\lambda')+n_{t+1}-v_{t+1}
				&=-\nabla_yG(\lambda)+n_{t+1}-v_{t+1}
				-\bigl(\nabla_yG(\lambda')-\nabla_yG(\lambda)\bigr).
			\end{align*}
			Hence $\|a+b\|^2-\|a\|^2
			\le\eta\|a\|^2+(1+\eta^{-1})\|b\|^2$ with
			$\eta=h\ell_\gamma/64$ for the first identity and
			$\eta=3hs_\gamma/8192$ for the second, together with
			\eqref{eq:q-lambda-gradient}, yields
			\begin{align*}
				&\frac1{\ell_\gamma}\|\nabla_xG(\lambda')+\xi_{t+1}\|^2
				-\frac1{\ell_\gamma}\|\nabla_xG(\lambda)+\xi_{t+1}\|^2\\
				&\le \frac h{64}\|\nabla_xG(\lambda)+\xi_{t+1}\|^2
				+\frac{J_c^2}{\ell_\gamma}
				\left(1+\frac{64}{h\ell_\gamma}\right)d^2,\notag\\
				&\frac1{\ell_\gamma}\|-\nabla_yG(\lambda')+n_{t+1}-v_{t+1}\|^2
				-\frac1{\ell_\gamma}\|-\nabla_yG(\lambda)+n_{t+1}-v_{t+1}\|^2\\
				&\le \frac{3hs_\gamma}{8192\ell_\gamma}
				\|-\nabla_yG(\lambda)+n_{t+1}-v_{t+1}\|^2+\frac{J_c^2}{\ell_\gamma}
				\left(1+\frac{8192}{3hs_\gamma}\right)d^2.
			\end{align*}
			Moreover, \eqref{eq:saddle-lambda-sensitivity} gives
			\[
			\frac{s_\gamma}{16\ell_\gamma}
			\bigl|\langle v_{t+1},y^\star(z,\lambda')-y^\star(z,\lambda)\rangle\bigr|
			\le\frac{s_\gamma J_c}{16\ell_\gamma\mu_y}\|v_{t+1}\|d.
			\]
			Young's inequality and $s_\gamma^2=2\ell_\gamma\mu_y$ give
			\[
			\frac{s_\gamma J_c}{16\ell_\gamma\mu_y}\|v_{t+1}\|d
			\le\frac{7hs_\gamma}{65536\ell_\gamma}\|v_{t+1}\|^2+\frac{128J_c^2}{7hs_\gamma\mu_y}d^2.
			\]
			Likewise, with
			$a=y_{t+1}-y^\star(z,\lambda)$ and
			$b=y^\star(z,\lambda')-y^\star(z,\lambda)$, we have
			\[
			\frac{\mu_y}{256}\bigl(\|a-b\|^2-\|a\|^2\bigr)
			\le\frac{h\mu_y s_\gamma}{2048}\|a\|^2
			+J_c^2\left(\frac1{32hs_\gamma\mu_y}
			+\frac1{256\mu_y}\right)d^2.
			\]
			Together with the value estimates above, this proves \eqref{eq:E-lambda-switch}.
			
			For \eqref{eq:B-lambda-switch}, apply
			$\|a+b\|^2\le\frac54\|a\|^2+5\|b\|^2$ in \eqref{eq:Bt}. The same perturbation bounds then give \eqref{eq:B-lambda-switch} directly from the definition of $K_B$ in \eqref{eq:KB}.
		\end{proof}
		
		\begin{theorem}\label{thm:solver}
			Suppose Assumption~\ref{ass:basic} holds, $\gamma\ge1$, $0\le\rho\le1$, and
			$\mu_y=\sigma_y+\rho\in(0,\ell_\gamma]$. If $\sigma_y>0$, also suppose Assumption~\ref{ass:strong-y} holds. Define the aggregate multiplier-drift coefficient
			\begin{equation}\label{eq:Kdr}
				K_{\rm dr}:=K_E+\frac12K_B
				+\frac{3\beta\ell_\gamma J_c^2}{16\mu_y^2}.
			\end{equation}
			Choose the multiplier relaxation as
			\begin{equation}\label{eq:theta-choice}
				\bar\theta_{\gamma,\mu_y}:=
				\frac{1}{2(1+C_c)\sqrt{1+K_{\rm dr}}},
				\qquad
				\theta_t:=\frac{\bar\theta_{\gamma,\mu_y}}{(t+1)^2}.
			\end{equation}
			Run Algorithm~\ref{alg:cpdeg} with
			\eqref{eq:aux-parameters},
			\eqref{eq:solver-parameters}, and
			\eqref{eq:theta-choice} and the initialization in Algorithm~\ref{alg:cpdeg}: $c(x_0,y_0)\le0$, $\lambda_0=0$, $z_0=x_0$, and $\xi_0=n_0=v_0=0$. Then $V_t\ge0$. With
			$d_t:=\|\lambda_{t+1}-\lambda_t\|$, one has
			\begin{equation}\label{eq:sumS}
				\sum_{t=0}^{T-1}S_t
				\le\frac{128\ell_\gamma}{\beta}
				\left(V_0+4C_c\sum_{t=0}^{T-1}d_t
				+K_{\rm dr}\sum_{t=0}^{T-1}d_t^2\right).
			\end{equation}
			Consequently, $S_j$ satisfies
			\begin{equation}\label{eq:solver-residual}
				\sqrt{S_j}\le \varepsilon_T:=
				\sqrt{\frac{128\ell_\gamma}{\beta T}
					\left(V_0+4C_c\sum_{t=0}^{T-1}d_t
					+K_{\rm dr}\sum_{t=0}^{T-1}d_t^2\right)}.
			\end{equation}
			For the choice \eqref{eq:theta-choice},
			\begin{equation}\label{eq:uniform-drift-budget}
				V_0+4C_c\sum_{t=0}^{T-1}d_t
				+K_{\rm dr}\sum_{t=0}^{T-1}d_t^2
				\le V_0+4C_c^2+1.
			\end{equation}
		\end{theorem}
		
		\begin{proof}{Proof.}
			\noindent
			Since $K_{\rm dr}\ge0$ and $\gamma\ge1$, \eqref{eq:theta-choice} gives $0<\theta_t\le1/2\le\gamma$ and $\sum_{t\ge0}\theta_t\le2\bar\theta_{\gamma,\mu_y}\le1$. Thus Lemma~\ref{lem:lambda-state} applies.
			For every fixed multiplier,
			Lemma~\ref{lem:energy-lower} gives
			\[
			\mathcal E_z^\lambda(w)+\frac{s_\gamma}{256\ell_\gamma^2}\|v\|^2\ge0.
			\]
			Compactness of $\mathcal X\times\mathcal Y$ and the multiplier bound in
			Lemma~\ref{lem:lambda-state} make $\inf_z p_\rho(z,\lambda_t)$ finite. Hence $V_t\ge0$.
			
			Apply \eqref{eq:fixed-center} with $\lambda_t$, then switch from
			$\lambda_t$ to $\lambda_{t+1}$ at the fixed center $z_t$ by
			\eqref{eq:E-lambda-switch}. Since $\nabla_zp_\rho(\cdot,\lambda_{t+1})$ is
			$6\ell_\gamma$-Lipschitz, we get
			\begin{align*}
				p_\rho(z_{t+1},\lambda_{t+1})-p_\rho(z_t,\lambda_{t+1})
				&\le \left\langle \nabla_zp_\rho(z_t,\lambda_{t+1}),z_{t+1}-z_t\right\rangle
				+3\ell_\gamma\|z_{t+1}-z_t\|^2.
			\end{align*}
			Substituting
			\[
			z_{t+1}-z_t
			=\beta\left(x_{t+1}-x^\star(z_t,\lambda_{t+1})
			-\frac{1}{2\ell_\gamma}\nabla_z p_\rho(z_t,\lambda_{t+1})\right)
			\]
			into this inequality, and using
			$\langle a,b\rangle\le \|a\|^2/(8\ell_\gamma)+2\ell_\gamma\|b\|^2$, yields
			\begin{align*}
				&p_\rho(z_{t+1},\lambda_{t+1})
				-p_\rho(z_t,\lambda_{t+1})\le-\frac{\beta}{4\ell_\gamma}
				\|\nabla_z p_\rho(z_t,\lambda_{t+1})\|^2
				+\frac{5\beta\ell_\gamma}{2}
				\|x_{t+1}-x^\star(z_t,\lambda_{t+1})\|^2.
			\end{align*}
			Adding \eqref{eq:center-motion} and using
			\[
			\frac{5\beta\ell_\gamma}{2}+\frac{1408\ell_\gamma\beta^2}{hs_\gamma}
			\le\frac{7h\ell_\gamma s_\gamma}{2048},
			\qquad
			\frac{\beta}{4\ell_\gamma}-\frac{352\beta^2}{\ell_\gamma hs_\gamma}
			\ge\frac{\beta}{8\ell_\gamma}
			\]
			yields
			\begin{align}
				&\left[p_\rho(z_{t+1},\lambda_{t+1})-
				\inf_z p_\rho(z,\lambda_{t+1})
				+\mathcal E_{z_{t+1}}^{\lambda_{t+1}}(w_{t+1})\right]\notag\\
				&-
				\left[p_\rho(z_t,\lambda_{t+1})-
				\inf_z p_\rho(z,\lambda_{t+1})
				+\mathcal E_{z_t}^{\lambda_{t+1}}(w_{t+1})\right]\notag\\
				&\le\frac12\mathcal B_t^{\lambda_{t+1}}
				-\frac{\beta}{8\ell_\gamma}
				\|\nabla_z p_\rho(z_t,\lambda_{t+1})\|^2.
				\label{eq:center-total}
			\end{align}
			Using \eqref{eq:B-lambda-switch} and \eqref{eq:B-le-D}, we obtain
			\begin{align}
				V_{t+1}-V_t
				&\le-\frac5{16}\mathcal D_t
				-\frac{\beta}{8\ell_\gamma}
				\|\nabla_z p_\rho(z_t,\lambda_{t+1})\|^2
				+4C_cd_t+\left(K_E+\frac12K_B\right)d_t^2.
				\label{eq:dynamic-descent}
			\end{align}
			Define
			\begin{align*}
				R_{x,t}
				&=A_{x,t}-2\ell_\gamma(x_{t+1}-x_t^\star)
				+\nabla_z p_\rho(z_t,\lambda_t),\quad
				R_{y,t}=A_{y,t}+v_{t+1}.
			\end{align*}
			Moreover, \eqref{eq:saddle-lambda-sensitivity} implies
			\[
			\|\nabla_z p_\rho(z_t,\lambda_t)\|^2
			\le2\|\nabla_z p_\rho(z_t,\lambda_{t+1})\|^2
			+\frac{8\ell_\gamma^2J_c^2}{\mu_y^2}d_t^2.
			\]
			Using $\|a+b+c\|^2\le3(\|a\|^2+\|b\|^2+\|c\|^2)$,
			$\|a+b\|^2\le2(\|a\|^2+\|b\|^2)$, we obtain 
			\begin{align}
				S_t
				&\le\frac{128\ell_\gamma}{\beta}
				\left(\frac5{16}\mathcal D_t
				+\frac{\beta}{8\ell_\gamma}
				\|\nabla_z p_\rho(z_t,\lambda_{t+1})\|^2\right)
				+\frac{24\ell_\gamma^2J_c^2}{\mu_y^2}d_t^2.\label{eq:S-from-D}
			\end{align}
			Summing \eqref{eq:dynamic-descent}, using $V_T\ge0$, and substituting into
			\eqref{eq:S-from-D} proves \eqref{eq:sumS}; the last term in
			\eqref{eq:Kdr} is exactly the rescaling of the final term in
			\eqref{eq:S-from-D}. The definition of $j$ proves
			\eqref{eq:solver-residual}.
			
			Finally, Lemma~\ref{lem:lambda-state} gives $d_t\le C_c\theta_t$. 
			$\sum_{t=0}^\infty\frac1{(t+1)^2}\le2$,
			$\sum_{t=0}^\infty\frac1{(t+1)^4}\le2$
			and \eqref{eq:theta-choice} imply
			\[
			\sum_td_t\le C_c,
			\qquad
			K_{\rm dr}\sum_td_t^2
			\le\frac{K_{\rm dr}C_c^2}
			{2(1+C_c)^2(1+K_{\rm dr})}\le1.
			\]
			This proves \eqref{eq:uniform-drift-budget}.
		\end{proof}

		Theorem~\ref{thm:solver} bounds the stationarity residual. Since $\xi_{j+1}$ and $n_{j+1}$ are normal vectors, removing the perturbation costs at most $\rho D_Y$:
		\begin{equation}\label{eq:q-residual}
			\begin{split}
				&\left[
				\dist^2\!\left(0,\nabla_x\mathcal L_\gamma(x_{\rm out},y_{\rm out},\lambda_j)
				+N_{\mathcal X}(x_{\rm out})\right)\right.\\[-2mm]
				&\qquad\left.+
				\dist^2\!\left(0,-\nabla_y\mathcal L_\gamma(x_{\rm out},y_{\rm out},\lambda_j)
				+N_{\mathcal Y}(y_{\rm out})\right)
				\right]^{1/2}
				\le \varepsilon_T+\rho D_Y.
			\end{split}
		\end{equation}
		
		\begin{theorem}[KKT guarantee]\label{thm:master}
			Suppose Assumption~\ref{ass:basic} holds, $0\le\rho\le1$, and Algorithm~\ref{alg:cpdeg} is run under Theorem~\ref{thm:solver}. If $\varepsilon_T\le1$, then its output with
			$\lambda_{\rm out}=\widetilde\lambda_{j+1}
			=[\lambda_j+\gamma c(x_{\rm out},y_{\rm out})]_+$ satisfies the following bounds 
			\begin{align}
				\max\{\mathcal R_x(x_{\rm out},y_{\rm out},\lambda_{\rm out}),
				\mathcal R_y(x_{\rm out},y_{\rm out},\lambda_{\rm out})\}
				&\le \varepsilon_T+\rho D_Y,\label{eq:master-station}\\
				\mathcal R_c(x_{\rm out},y_{\rm out})
				&\le\frac{M_{\rm aug}}{\gamma},\label{eq:master-feas}\\
				\mathcal R_{\rm comp}(x_{\rm out},y_{\rm out},\lambda_{\rm out})
				&\le\frac{M_{\rm aug}(M_{\rm aug}+C_c)}{\gamma}.
				\label{eq:master-comp}
			\end{align}
			where $M_{\rm aug}$ defined in \eqref{eq:Mmu} is the uniform upper bound on $\|\lambda_{\rm out}\|_1$.
			Consequently, all four residuals are at most $\varepsilon$ whenever
			\begin{equation}\label{eq:master-choice}
				\varepsilon_T+\rho D_Y\le\varepsilon,
				\qquad
				\gamma\ge
				\frac{\max\{M_{\rm aug},M_{\rm aug}(M_{\rm aug}+C_c)\}}
				{\varepsilon}.
			\end{equation}
		\end{theorem}
		
		\begin{proof}{Proof.}
			\noindent
			By Lemma~\ref{lem:normals} and \eqref{eq:St}, we have
			\[
			\begin{aligned}
				\sqrt{S_j}
				&=\left[
				\|\nabla_x\mathcal L_\gamma(x_{\rm out},y_{\rm out},\lambda_j)+\xi_{j+1}\|^2\right.\\
				&\quad\left.+\|-\nabla_y\mathcal L_\gamma(x_{\rm out},y_{\rm out},\lambda_j)
				+\rho y_{\rm out}+n_{j+1}\|^2
				\right]^{1/2}
				\le\varepsilon_T.
			\end{aligned}
			\]
			Thus, $\rho y_{\rm out}$ is the perturbation in the $y$-stationarity residual. Substituting \eqref{eq:gradqx}--\eqref{eq:gradqy} and
			$\lambda_{\rm out}=\widetilde\lambda_{j+1}$, and using $\|y_{\rm out}\|\le D_Y$, gives
			\[
			\max\{\mathcal R_x(x_{\rm out},y_{\rm out},\lambda_{\rm out}),
			\mathcal R_y(x_{\rm out},y_{\rm out},\lambda_{\rm out})\}
			\le \varepsilon_T+\rho D_Y,
			\]
			which is \eqref{eq:master-station}. Since $\varepsilon_T\le1$, Lemmas~\ref{lem:lambda-state} and~\ref{lem:mu-bound}, applied with $\lambda=\lambda_j$, $\mu=\lambda_{\rm out}$, and $r=\varepsilon_T$, yield
			$\|\lambda_{\rm out}\|_1\le M_{\rm aug}$. Substituting this bound into Lemma~\ref{lem:feas-comp} gives
			\[
			\mathcal R_c(x_{\rm out},y_{\rm out})
			\le\frac{M_{\rm aug}}{\gamma},
			\qquad
			\mathcal R_{\rm comp}(x_{\rm out},y_{\rm out},\lambda_{\rm out})
			\le\frac{M_{\rm aug}(M_{\rm aug}+C_c)}{\gamma},
			\]
			which are \eqref{eq:master-feas} and \eqref{eq:master-comp}. Finally, choose $\rho$ and $T$ so that $\varepsilon_T+\rho D_Y\le\varepsilon$, and choose
			\[
			\gamma\ge
			\frac{\max\{M_{\rm aug},M_{\rm aug}(M_{\rm aug}+C_c)\}}{\varepsilon}.
			\]
			Then the right-hand sides of \eqref{eq:master-station}--\eqref{eq:master-comp} are at most $\varepsilon$, proving \eqref{eq:master-choice} and the result.
		\end{proof}

		\subsection{Nonconvex-(Strongly) Concave Complexity}\label{sec:sc}
		We now derive iteration bounds for the strongly concave and concave cases.
	A complete first-order oracle evaluation at $(x,y)$ returns
	$f(x,y)$, $c(x,y)$, $\nabla_x f(x,y)$, $\nabla_y f(x,y)$,
	$J_xc(x,y)$, and $J_yc(x,y)$.
	By caching these quantities, the surrogate gradients can be
	recomputed after multiplier and center updates without additional
	oracle calls. Thus, one initial evaluation and two evaluations
	per iteration, at the prediction and correction points, suffice,
	giving a total of $1+2T$ evaluations for $T$ iterations.
	This count excludes projections, algebraic operations, and
	the construction of the required feasible initial point.

		\begin{theorem}[$O(\varepsilon^{-7/2})$ complexity under strong concavity]
			\label{thm:sc}
			Suppose Assumptions~\ref{ass:basic} and~\ref{ass:strong-y} hold.
			Set $C_{\rm aug}:=\max\{1,M_{\rm aug},
			M_{\rm aug}(M_{\rm aug}+C_c)\}$.
			For $\varepsilon\in(0,1]$, choose
			\begin{equation}\label{eq:strong-choice}
				\gamma_\varepsilon:=\frac{C_{\rm aug}}{\varepsilon},
				\quad \rho_\varepsilon:=0,\quad
				\ell_\varepsilon:=\max\left\{1,
				L_f+L_{Jc}(C_c+\gamma_\varepsilon C_c)
				+\gamma_\varepsilon J_cL_c,\sigma_y\right\}.
			\end{equation}
			Initialize with any $x_0\in\mathcal X$ and $y_0\in\mathcal Y$
			satisfying $c(x_0,y_0)\le0$, and set the remaining parameters by
			\eqref{eq:aux-parameters}, \eqref{eq:solver-parameters}, and
			\eqref{eq:theta-choice}, with $\mu_y=\sigma_y$.
			Set $\overline B_{\sigma_y}:=B_{\sigma_y}+4C_c^2+1$ and
			\begin{equation}\label{eq:Tstrong}
				T:=\max\left\{1,
				\left\lceil
				\frac{512\ell_\varepsilon\overline B_{\sigma_y}}
				{\beta\varepsilon^2}
				\right\rceil\right\}.
			\end{equation}
		Then Algorithm~\ref{alg:cpdeg} returns an
		$\varepsilon$-KKT point using $1+2T$
		complete first-order evaluations. This evaluation count is of order $\mathcal O(\varepsilon^{-7/2})$.
		\end{theorem}
		
	\begin{proof}{Proof.}
		Lemma~\ref{lem:initial-bound}, Theorem~\ref{thm:solver},
		and \eqref{eq:Tstrong} give
		$\varepsilon_T\le\varepsilon/2\le1$.
		Since $\rho_\varepsilon=0$, Theorem~\ref{thm:master}
		bounds the stationarity residual by $\varepsilon/2$.
		The definition of $C_{\rm aug}$ and the choice
		$\gamma_\varepsilon=C_{\rm aug}/\varepsilon$
		bound feasibility and complementarity by $\varepsilon$.
		Thus the returned point is an $\varepsilon$-KKT point.
		For fixed problem data, the parameter choices imply that
		$\beta$ scales as $\sqrt{\sigma_y/\ell_\varepsilon}$,
		while $\overline B_{\sigma_y}$ is independent of $\varepsilon$.
		Hence \eqref{eq:Tstrong} bounds the evaluation count $1+2T$
		by order
		$\mathcal O(1+\ell_\varepsilon^{3/2}
		\overline B_{\sigma_y}/(\sqrt{\sigma_y}\,\varepsilon^2))$.
		Since \eqref{eq:strong-choice} gives
		$\ell_\varepsilon$ of order $\mathcal O(\varepsilon^{-1})$,
		the claimed $\mathcal O(\varepsilon^{-7/2})$ bound follows.
	\end{proof}

\begin{theorem}[$O(\varepsilon^{-4})$ complexity under concavity]
	\label{thm:concave}
	Suppose Assumption~\ref{ass:basic} holds and $f(x,\cdot)$ is concave.
	Set $C_{\rm aug}:=\max\{1,M_{\rm aug},
	M_{\rm aug}(M_{\rm aug}+C_c)\}$ and $D_+:=\max\{D_Y,1\}$.
	For $\varepsilon\in(0,1]$, choose
	\begin{equation}\label{eq:concave-choice}
		\gamma_\varepsilon:=\frac{C_{\rm aug}}{\varepsilon},
		\quad \rho_\varepsilon:=\frac{\varepsilon}{2D_+},
		\quad \ell_\varepsilon:=\max\left\{1,
		L_f+L_{Jc}(C_c+\gamma_\varepsilon C_c)
		+\gamma_\varepsilon J_cL_c\right\}.
	\end{equation}
	Initialize with any feasible pair $(x_0,y_0)\in\mathcal X\times\mathcal Y$,
	and set the remaining parameters by
	\eqref{eq:aux-parameters},
	\eqref{eq:solver-parameters}, and
	\eqref{eq:theta-choice}, with $\sigma_y=0$ and $\mu_y=\rho_\varepsilon$.
	Set $\overline B_0:=B_0+4C_c^2+1$ and
	\begin{equation}\label{eq:Tconcave}
		T:=\max\left\{1,
		\left\lceil
		\frac{512\ell_\varepsilon\overline B_0}
		{\beta\varepsilon^2}
		\right\rceil\right\}.
	\end{equation}
	Then Algorithm~\ref{alg:cpdeg} returns an
	$\varepsilon$-KKT point using $1+2T$
	complete first-order evaluations. This evaluation count is of order $	\mathcal O\!\left(
	1+\frac{\ell_\varepsilon^{3/2}\overline B_0}
	{\sqrt{\rho_\varepsilon}\,\varepsilon^2}\right),$
	and hence of order $\mathcal O(\varepsilon^{-4})$.
\end{theorem}
\begin{proof}{Proof.}
	Lemma~\ref{lem:initial-bound}, \eqref{eq:solver-residual},
	and \eqref{eq:Tconcave} give
	$\varepsilon_T\le\varepsilon/2\le1$.
	Since $\rho_\varepsilon D_Y\le\varepsilon/2$,
	Theorem~\ref{thm:master} bounds both stationarity residuals
	by $\varepsilon$.
	The definition of $C_{\rm aug}$ and the choice
	$\gamma_\varepsilon=C_{\rm aug}/\varepsilon$
	bound feasibility and complementarity by $\varepsilon$.
	Thus the returned point is an $\varepsilon$-KKT point.
	For fixed problem data, the parameter choices imply that
	$\beta$ scales as $\sqrt{\rho_\varepsilon/\ell_\varepsilon}$,
	while $\overline B_0$ is independent of $\varepsilon$.
	Hence \eqref{eq:Tconcave} bounds the evaluation count $1+2T$
	by order $\mathcal O(1+\ell_\varepsilon^{3/2}\overline B_0/
	(\sqrt{\rho_\varepsilon}\,\varepsilon^2))$.
	Since \eqref{eq:concave-choice} gives
	$\ell_\varepsilon$ of order $\mathcal O(\varepsilon^{-1})$
	and $\rho_\varepsilon=\varepsilon/(2D_+)$,
	the claimed $\mathcal O(\varepsilon^{-4})$ bound follows.
\end{proof}

		\section{Numerical Results}\label{sec:nu}
		
		We compare ALDER with MGD~\cite{Tsaknakis}, PGmsAD~\cite{Dai},
		and FAL~\cite{lu23} on a nonlinearly constrained synthetic problem
		and a fairness-aware classification problem. Performance is measured
		by the current KKT residual $\max\{\mathcal R_x,\mathcal R_y,\mathcal R_c,\mathcal R_{\rm comp}\}$ in Definition~\ref{def:kkt-gs}
		versus solver CPU time. These experiments assess the computational
		efficiency of the methods at prescribed residual thresholds.
		All experiments are implemented in Python~3.14.0
		and run on an Apple M4 Pro computer with 24~GB of RAM.
		
	\subsection{Synthetic Example}
	
	We consider Example~1.1 from~\cite{CaoZengZhangZhou2026},
	\begin{equation}\label{eq:numerical-linear}
		\min_{x\in\mathcal X}
		\max_{\substack{y\in\mathcal Y\\e^\top y-\|x\|^2\le0}}
		\left(\frac{\|x\|^2}{2}-1\right)^2
		+2\|x\|^2-(e^\top x)^2
		-\frac{\|y-e\|^2}{2}+\frac{x^\top y}{2},
	\end{equation}
	where $\mathcal X=[-5/4,5/4]^2$,
	$\mathcal Y=[-10,10]^2$, and $e=(1,1)^\top$.
	
	All four methods start from a common feasible initialization.
	For ALDER, we use
	$(\gamma,\allowbreak h,\allowbreak\alpha,\allowbreak
	\vartheta,\allowbreak\beta)
	=(1.005,\allowbreak0.16,\allowbreak0.25,\allowbreak
	3.5,\allowbreak0.5)$
	and $\theta_t=\min\{\gamma,10^8/(t+1)^2\}$.
	For MGD, following the experimental setup
	in~\cite{CaoZengZhangZhou2026}, we set
	$(x\text{-step},\allowbreak y\text{-step},\allowbreak
	\text{dual step},\allowbreak\text{inner steps})
	=(0.01,\allowbreak0.01,\allowbreak0.1,\allowbreak1)$.
	For PGmsAD, the corresponding parameters are
	$(0.3,\allowbreak0.5,\allowbreak0.3,\allowbreak2)$.
	For FAL, we use
	$(\tau,\allowbreak\text{smoothness scale},\allowbreak
	\text{diameter scale})
	=(0.999,\allowbreak0.02,\allowbreak0.1)$.
	The parameters of ALDER, PGmsAD, and FAL are selected
	using equal CPU-time budgets for parameter tuning. The practical parameter choices used here need not satisfy
	the sufficient conditions of the complexity analysis.

		\begin{figure}[!ht]
			\centering
			\includegraphics[width=0.65\textwidth]{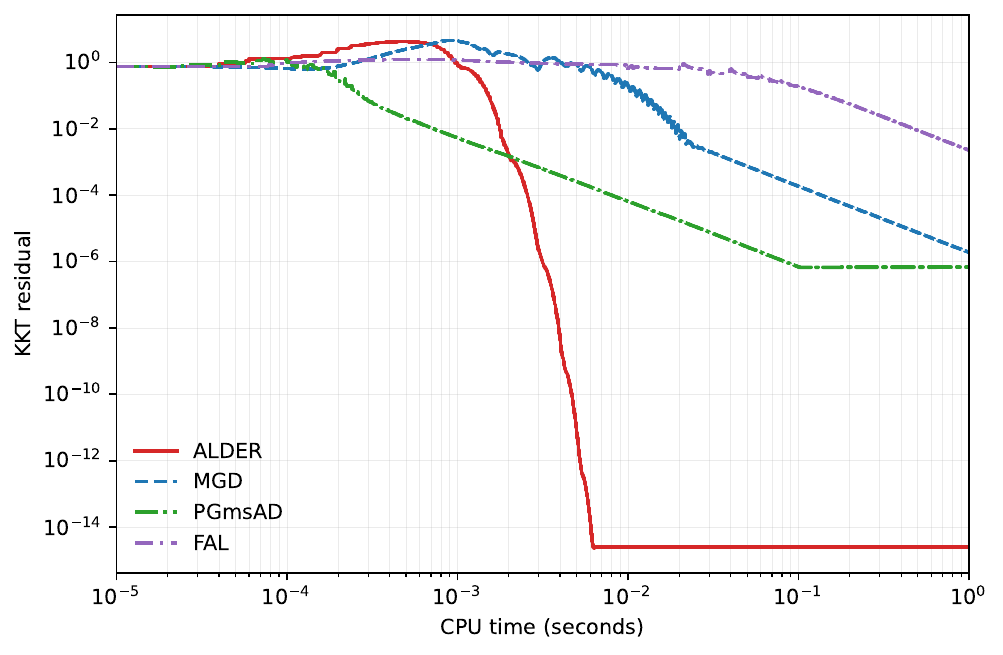}
			\caption{KKT residual versus solver CPU time for
			 problem~\eqref{eq:numerical-linear}.}
			\label{fig:linear-coupled-cpu}
		\end{figure}
		
		As shown in Figure~\ref{fig:linear-coupled-cpu}, ALDER reaches
		a KKT residual of $10^{-4}$ in the shortest CPU time among
		the four methods. This result highlights its computational
		efficiency at the tested accuracy on this instance.
		
		\subsection{Fairness-aware Classification}
	We consider fairness-aware classification on the UCI Adult
	data set~\cite{Adult}, following the adversarial construction
	of Zhang et al.~\cite{ZhangLemoineMitchell2018} and the constrained
	model of Cao et al.~\cite{CaoZengZhangZhou2026},
	\begin{equation}\label{eq:adult-fairness}
		\min_{\theta\in[-3,3]^d}
		\max_{\substack{\phi\in[-5,5]^3\\
				L_{\rm adv}(\theta,\phi)\le\kappa}}
		L_{\rm cls}(\theta)-\nu L_{\rm adv}(\theta,\phi).
	\end{equation}
	The classification and adversarial losses are
	\[
	L_{\rm cls}(\theta)
	=\frac1N\sum_{i=1}^N\ell(g_\theta(a_i),b_i),
	\qquad
	L_{\rm adv}(\theta,\phi)
	=\frac1N\sum_{i=1}^N
	\ell\!\left(\phi^\top(p_\theta(a_i),b_i,1)^\top,s_i\right),
	\]
	where $a_i\in\mathbb R^{102}$, $b_i\in\{0,1\}$, and
	$s_i\in\{0,1\}$ denote the feature vector, class label, and
	sensitive attribute of sample $i$, respectively, and
	$\ell(q,r)=\log(1+e^q)-rq$ is the logistic loss.
	The predictor is defined by
	\[
	g_\theta(a_i)=w^\top\tanh(W^\top a_i+b)+b_0,
	\qquad
	p_\theta(a_i)=\frac{1}{1+e^{-g_\theta(a_i)}},
	\]
	where $W\in\mathbb R^{102\times4}$, $w,b\in\mathbb R^4$,
	$b_0\in\mathbb R$, and $\tanh$ is applied componentwise.
	The parameter vector $\theta$ collects all predictor parameters,
	so $d=417$.
	We set $\nu=0.5$ and $\kappa=0.69$.
	
\begin{figure}[!ht]
			\centering
			\includegraphics[width=0.65\textwidth]{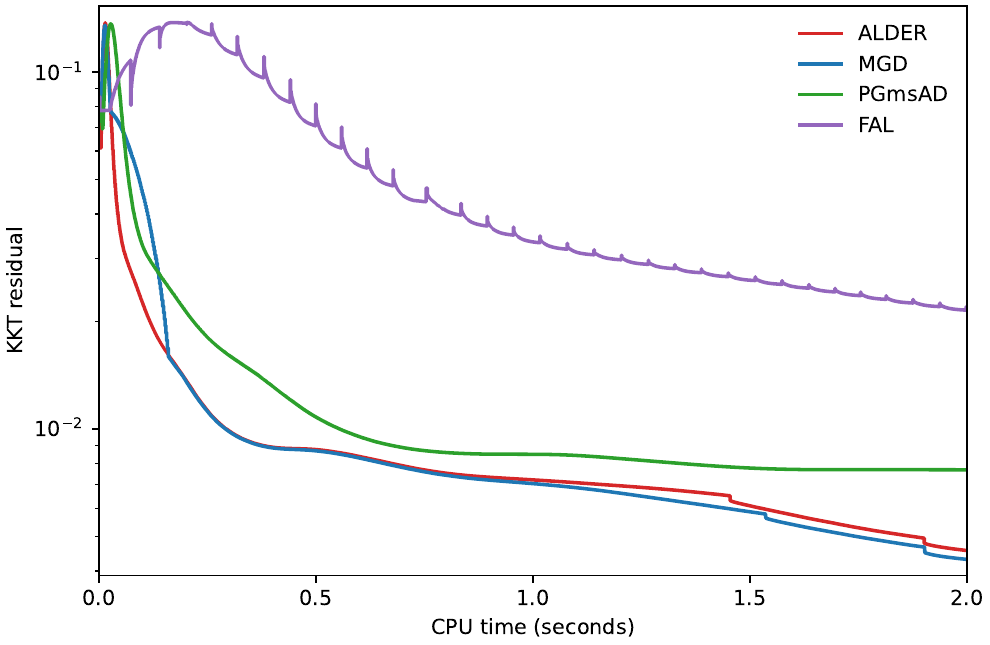}
			\caption{KKT residual versus solver CPU time for
				problem~\eqref{eq:adult-fairness}.}
			\label{fig:adult-fairness}
		\end{figure}

	All four methods start from a common feasible initialization,
	and their parameters are selected using equal CPU-time budgets
	for parameter tuning.
	For ALDER, we use
	$(\gamma,\allowbreak\rho,\allowbreak h,\allowbreak
	\alpha,\allowbreak\vartheta,\allowbreak\beta)
	=(1.005,\allowbreak10^{-4},\allowbreak0.4,\allowbreak
	0.8,\allowbreak0.2,\allowbreak0.2)$
	and $\theta_t=\min\{\gamma,1000/(t+1)^{1.5}\}$.
	For MGD, we set
	$(x\text{-step},\allowbreak y\text{-step},\allowbreak
	\text{dual step},\allowbreak\text{inner steps})
	=(0.2,\allowbreak0.15,\allowbreak0.1,\allowbreak5)$.
	For PGmsAD, the corresponding parameters are
	$(0.2,\allowbreak0.15,\allowbreak50.2,\allowbreak1)$.
	For FAL, we use $\tau=0.995$, with smoothness and diameter
	scales of $0.005$ and $0.2$, respectively.
Figure~\ref{fig:adult-fairness} shows that ALDER and MGD exhibit
comparable performance, both reaching the tested KKT residual
threshold in less CPU time than PGmsAD and FAL.
		
		\section{Conclusions}\label{section_conclu}
		
		We developed ALDER, a single-loop augmented-Lagrangian damped
		extragradient method for smooth nonconvex--concave minimax problems
		with nonlinear coupled inequality constraints. The method combines
		projected prediction--correction steps with an explicit
		projected-and-relaxed multiplier update, avoiding inner subproblem
		solves and penalty-parameter updates. We established uniform
		boundedness of the relaxed multiplier sequence and, under extended
		MFCQ, boundedness of the projected output multiplier.
		With a feasible initialization and a zero initial multiplier,
		ALDER obtains an $\varepsilon$-KKT point using
		$\mathcal O(\varepsilon^{-7/2})$ complete first-order evaluations
		under strong concavity and $\mathcal O(\varepsilon^{-4})$
		under concavity. To the best of our knowledge, these bounds have
		the best known dependence on $\varepsilon$ among single- and
		multi-loop first-order methods for this problem class.
		In particular, they remove the logarithmic factor in the corresponding
		FAL bounds~\cite{LuMei2026SC,lu23} while retaining a single-loop
		structure.

	\end{document}